\documentclass[lettersize,journal]{IEEEtran}
\usepackage{amsmath,amsfonts}
\usepackage{algorithmic}
\usepackage{algorithm}
\usepackage{array}
\usepackage[caption=false,font=normalsize,labelfont=sf,textfont=sf]{subfig}
\usepackage{textcomp}
\usepackage{stfloats}
\usepackage{url}
\usepackage{verbatim}
\usepackage{graphicx}
\usepackage{xcolor}
\usepackage{cite}
\def\be{\begin{eqnarray}}
\def\ee{\end{eqnarray}}
\def\bee{\begin{eqnarray*}}
\def\eee{\end{eqnarray*}}
   
\newcommand{\ba}{\begin{align}} 
\newcommand{\ea}{\end{align}}

\usepackage{amsfonts} 
\usepackage{amssymb}

\def\be{\begin{eqnarray}}
\def\ee{\end{eqnarray}}
\def\bee{\begin{eqnarray*}}
\def\eee{\end{eqnarray*}}

\newtheorem{thm}{Theorem}

\newtheorem{lemma}[thm]{Lemma}

\def\wh{\widehat}

\def\E{\mathbb{E}}
\def\P{\mathbb{P}}

\def\bmx{\begin{pmatrix}}
\def\emx{\end{pmatrix}}
\def\bmat{\begin{pmatrix}}
\def\emat{\end{pmatrix}}

\def\mid{{\rm mid}}

\newcommand{\ignore}[1]{}

\begin{document}
\title{Personalised Feedback Control, Social Contracts, and Compliance Strategies for Ensembles}

\author{P.~Ferraro,\thanks{Pietro Ferraro, Lianna Zhao and Robert Shorten are with the Dyson School of Design Engineering, Imperial College London.} L.~Zhao,
  C.~King and R.~Shorten\thanks{Christopher King is with the Department of Mathematics, Northeastern University, Boston, MA 02115 USA. } }%

\maketitle

\begin{abstract}
This paper describes the use of Distributed Ledger Technologies as a means to create personalised social nudges and to influence the behaviour of agents in a smart city environment. Specifically, we present a scheme to price personalised risk in sharing economy applications. We provide proofs for the convergence of the proposed stochastic system and we validate our approach through the use of extensive Monte Carlo simulations. 
\end{abstract}

\begin{IEEEkeywords}
Social impacts, Smart Cities, Cyber-Physical Systems, Control Theory, Distributed Ledger Technologies.
\end{IEEEkeywords}

\section{Introductory remarks}
\label{Sec: Introduction}
{\color{black}

Our objective in this paper is to develop a framework to underpin the design of new forms of social contracts to guide the interaction of citizens, IoT-devices, and IoT enabled urban infrastructures. Specifically, we wish to develop control theoretic algorithms that can be deployed using 
distributed ledger technologies (DLTs) that nudge agents to respect social contracts in a range of Sharing Economy applications.  Such technologies are of great relevance in a range of Sharing Economy applications that often require the good behaviour and {\em compliance} of humans. For example, people should return shared assets in a timely manner as promised, and in good condition. Our basic proposal is to deploy a number of digital tokens as a bond or deposit to incentivise {\em compliance}: if the agent remains in compliance then these tokens are returned, otherwise they, or some part of them, are lost \cite{accesspaper}. As such, our algorithms can be viewed as mechanisms to realise personalised nudges that use financial rewards and penalties in order to direct how individual agents behave. Our contribution, in this specific context, is to develop DLT enabled-feedback strategies that manage the level of compliance based on personalised (but anonymous) interventions.\newline

\textcolor{black}{To provide some context on when the need for such a system may arise
we describe briefly a V2G (vehicle to grid) battery swapping architecture for shared e-bikes
that we have designed in our lab. Roughly speaking, our smart battery system is a unit that both charges 
and aggregates the batteries from a number of e-bikes. For example, this
could be part of a system in an apartment block that provides backup power
in the building for essential services, as well as charging individual batteries. 
To see how the need for our compliance system 
arises in the context of this example, we first note that the batteries are financially valuable
parts of the e-bike system. In our system, apartment owners would give residents
access to an e-bike. Residents purchase
tokens (whose value would exceed that of a battery)
and use a digital deposit system as described above; namely
in order to release a battery, users would deposit a token
into the ChargeWall, which would then be returned when the
discharged battery is returned (the social contract). In the remainder of this paper we describe methods to
set the correct number of tokens to encourage bike users to comply with the 
social contract of returning a battery to the ChargeWall. 
}

\begin{figure}[H]
\centering
\includegraphics[width=\columnwidth]{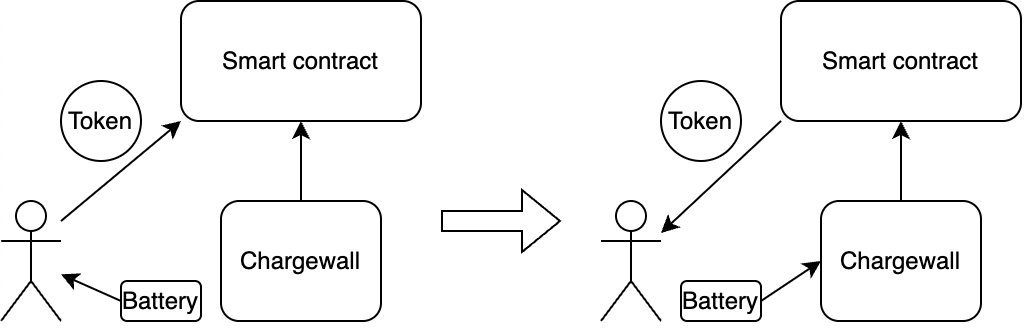}
\caption{Sequence to borrow a battery from ChargeWall. The user will deposit a certain amount of tokens, through a smart contract and these will be returned in full, once the battery is returned.}
\label{fig: Token}
\end{figure}

The issue of designing engineering systems to enforce compliance with social contracts is not new, and has indeed become very topical recently in a number of fields. In particular, in the context of Covid 19 (social distancing and mask wearing), several systems based on machine learning have been proposed to encourage compliance with social contracts. Much of this recent work on this topic involves using AI techniques to identify non-compliant actors, and to rely on some complementary strategy to enforce the social contract \cite{ML1}-\cite{ML4}. Our approach is quite different and does not separate identification of non-compliance from enforcement of the social contract. Rather, as we have mentioned, we use a personalised pricing strategy to encourage agents to comply with social contracts using financial rewards and penalties. Thus, our work is related to the general area of   \emph{Transactive Control}; namely, the use of financial transactions as a feedback signal to improve quality of service in various domains: some examples of work in this area can be found in \cite{Phan}-\cite{Junjie}. Other specific instances of dynamic pricing are \cite{Li} (incentivizing users to schedule electricity-consuming applications more prudently), \cite{Chekired} (managing  EVs charging and discharging in order to reduce the peak loads), \cite{Bejestani} (combining the classical hierarchical control in the power grid with market transactions) and \cite{Hao} (where the authors propose a transactive control system of commercial building heating, ventilation, and air-conditioning for demand response).\newline 

As we have mentioned, there are strong similarities between our work and transactive control. While this is true, our contribution differs from this prior work in this area in a number of ways. First, our proposal is a form of \emph{dynamic deposit pricing}: we explicitly price the {\em risk of not using an asset correctly}, rather than pricing the  {\em cost of access to the asset}. Second, the overall feedback strategy is based on a form of personalised nudging of agents that is akin to {\em polluter pays} type pricing models. In this specific context the contributions of this paper are to: (i) present a control system for stochastic agents, enabled by the use of DLTs
(Digital Ledger Technologies), to achieve a desired level of compliance in social contracts; and (ii) present a theoretical analysis that establishes the stochastic convergence of the proposed control system.}\newline

\subsection{Paper structure} The remainder of this paper is organised as follows: In Section II we discuss the proposed control architecture and why the use of Distributed Ledger Technologies is desirable when designing a compliance scheme. Sections III and IV describe the proposed control strategy and provide theoretical guarantees on its convergence. In Section V we show the effectiveness of the proposed approach through extended Monte Carlo simulations and, finally, Section VI summarizes the presented results and outlines future lines of research.\newline

\section{Control Architecture}
{\color{black}The basic idea presented in this paper is to design algorithms which implement social policies, using distributed ledger technology (often referred to as DLTs; one example is Blockchain). We combine the notions of digital identity and smart contracts coming from DLTs with the rigorous design methods afforded by control theory, in order to design personalised interventions for ensembles of agents.\newline  

In our context a {\bf Social Contract} is a set of rules, or policies, designed to govern the interaction between humans, other humans, and societal infrastructures. For example, in the context of shared assets such as a pool of vehicles, a social contract might require that vehicles are returned to a specified location at a contracted time. Another example of a social contract is the requirement that plastic bottles are returned to point-of-sale after use. \newline

The basic idea is to use digital tokens pegged at a stable fixed value, in the form of a cryptocurrency, to nudge users to comply with a social contract.
To be more specific, let $E$ represent a general statement, such as \emph{people should wear a face mask} or \emph{a shared vehicle needs to be returned at a specified time and place}. The agent purchases some number of digital tokens in order to participate in the social contract $E$: if the agent behaves in compliance with the contract then tokens are returned in their entirety, whereas if the agent does not then some tokens are lost. Thus, the risk of losing tokens is the mechanism that encourages agents to comply with these social contracts.\newline

A DLT is nothing more than a shared database. By its very nature it is decentralised and therefore no central authority is required in order to achieve consensus amongst users. In DLTs, transactions are pseudo-anonymous\protect\footnote{https://laurencetennant.com/papers/anonymity-iota.pdf}, and their content can be encrypted. This allows every agent to own and manage access to the data present in their own transactions. In our setting the only requirement is that the ownership of the tokens, used in the control, needs to remain visible to the compliance control algorithms, whereas other information (e.g., user quality of service, statistics on the usage of the system) can be encrypted. This allows each user to maintain ownership of their data and to use them as they please (e.g., to monetize them at a later stage). Token balances, and records of compliance on the ledger, associated with digital identities are the basis of the compliance system.  Finally, in order to enable the compliance-based control, we will focus on the use of DLTs built around Directed Acyclic Graphs (DAGs) such as IOTA\cite{Popov}\cite{TangleNew}\cite{Tangle}. These kinds of ledgers faciliate large transaction speeds, and are fee-less. In contrast, many standard payment systems (e.g., VISA, Mastercard) and classical Blockchain architectures (e.g., Bitcoin, Ethereum) require users to pay a fee for each transaction. This makes such systems inadequate to serve as the backbone for the proposed compliance scheme. It is essential that the deposited token is returned in its entirety to the owner in the event of full compliance with rule $E$: the proposed social compliance mechanism would break down if the agent were required to pay a fee every time a token is deposited or returned, as this would effectively erode the value of her tokens over time.}\newline

\begin{figure}[H]
\includegraphics[width=0.95\columnwidth]{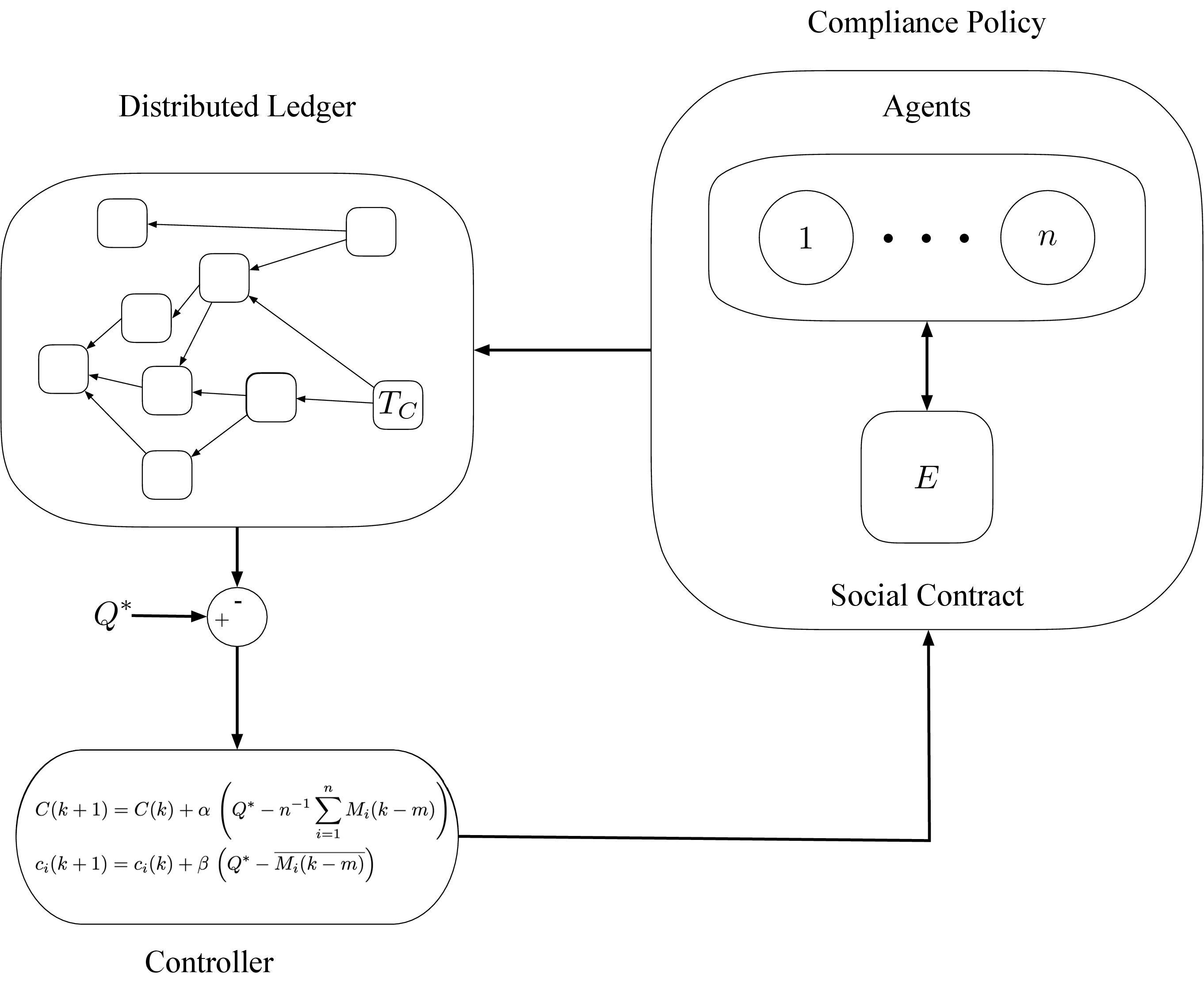}
\caption{The three components of the proposed compliance architecture are: the Distributed Ledger that acts as the communication backbone of the infrastructure, the compliance policy and the feedback controller.}
\label{fig: architecture}
\end{figure}

{\color{black}The proposed architecture is shown in Figure  \ref{fig: architecture}. The scheme is divided into three main components.\newline
\begin{itemize}
\item The Distributed Ledger discussed in this section, whose purpose is threefold: firstly, it acts as the communication backbone for the whole infrastructure, secondly it enables the deposit mechanism for the digital bond, thirdly it provides the controller with the current state of the network.\newline
\item The Physical Layer, in which agents interact with their environment in the setting of the social contract $E$.\newline
\item The Controller Layer, whose task is to regulate the price of the token bonds in order to achieve the desired level of compliance. Notice that the controller is not centralised as, due to the nature of the Distributed Ledger, every agent can implement the compliance scheme locally.
\end{itemize}}

 The latter two components of the architecture will be the focus of the next two Sections. A complete discussion of DAG-based DLTs and their comparison to classical Blockchain is beyond the scope of this paper; the interested reader can refer to \cite{accesspaper}  \cite{Popov} \cite{IoTPaper1} \cite{IoTPaper2}  \cite{TaCPaper} for a thorough discussion of their properties. For the purpose of this paper, all we need to assume is that there is a fast and secure way to execute the deposit and retrieval of these token bonds. \newline

\emph{Remark: }Before proceeding to the analysis and the modelling of the proposed framework it is worth stressing that the issue of compliance is often not incorporated into algorithms which are designed to regulate, control, and optimise city infrastructures. Many studies addressing human behaviour in this context assume full compliance with policies that have been engineered to optimally organise city infrastructures. As an example, consider traffic flow optimization: a crucial element that is often left out is that humans break rules, and the effect of this rule-breaking profoundly affects how cities operates and how well the engineered algorithms actually perform.

\section{Policy choices for compliance with social contracts}
\label{sec: Policy}
As explained before we are interested in using DLT's to create a type of {\em digital bond} which will encourage compliance with social contracts. Figure \ref{fig: Token}  provides a visual representation of this basic idea. In order to engage with this social scheme each agent stakes an amount of tokens (to which some monetary cost is associated) that acts as a bond. {\color{black}This is shown in Figure \ref{fig: Token}, by the arrow that goes from the agent to the transaction $T_C$ in the distributed ledger. The transaction $T_C$  represents the deposit of the tokens from the wallet of the agent to the wallet of a smart contract (these are distributed computer programs which execute as soon as certain conditions are met \cite{ethereum1}\cite{ethereum2}). Once the smart contract verifies that the agent complied with the rules of of social contract $E$, it returns the funds to the agent wallet (through a subsequent transaction). Notice that due to the nature of smart contracts,  the operations of deposit and return of the tokens are carried out automatically}. All these operations are recorded on a DLT that is shared amongst all agents (anonymously). \newline

\begin{figure}[H]
\centering
\includegraphics[width=0.75\columnwidth]{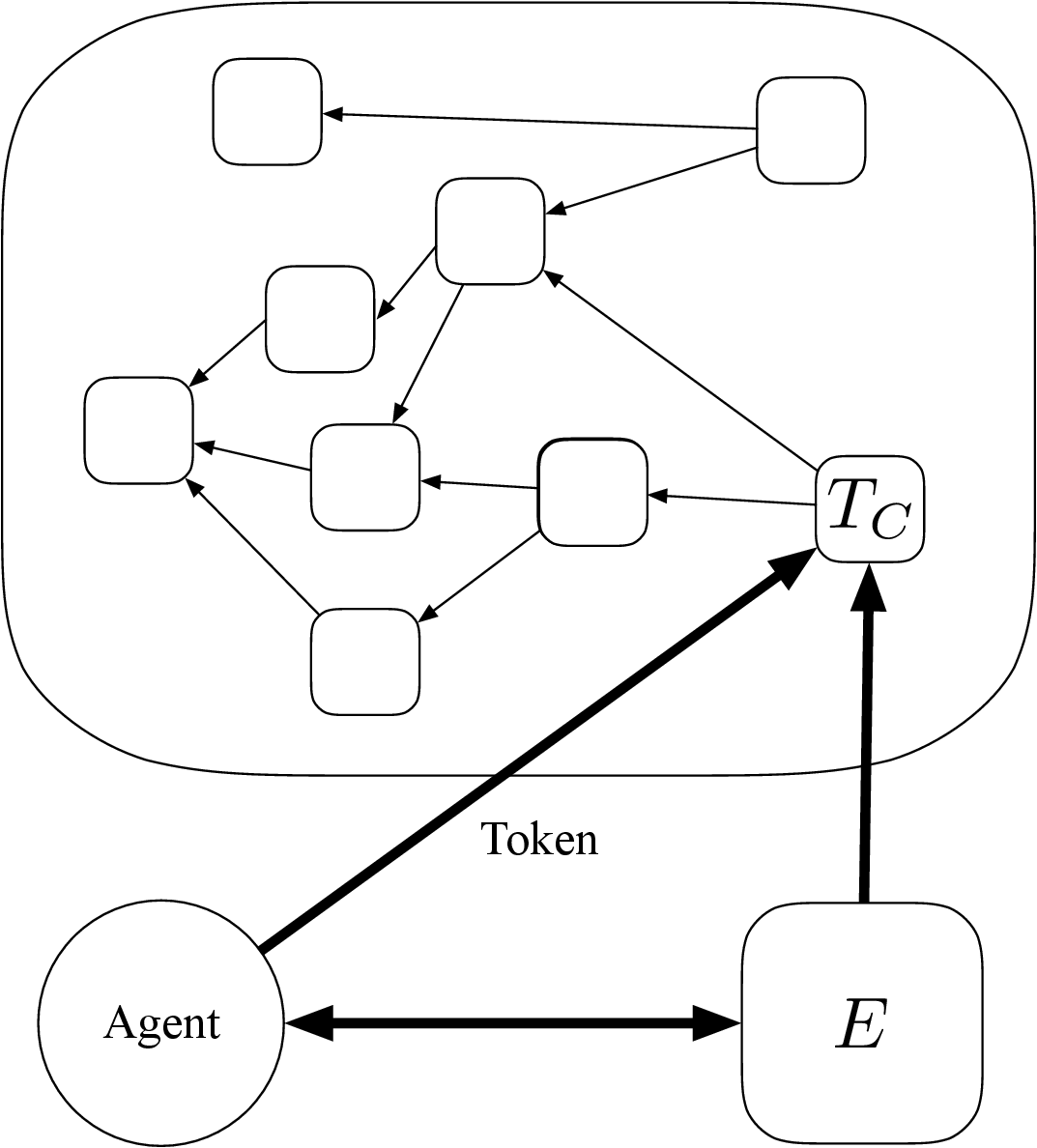}
\caption{Compliance Policy: each agent deposits a certain amount of tokens that are returned at the end of the activity. $T_C$ is a transaction recorded in a DAG-based DLT that certifies that the tokens have been deposited. Similarly if the agent complies, another transaction is recorded to return the tokens.}
\label{fig: Token}
\end{figure}

A basic question that arises is how to price this bond: namely, {\em how many tokens should be required as a bond in order to assure compliance with a social contract?} Clearly, if this number is too low, one can expect low levels of compliance, and if it is too high, activity will cease and the social contract will be meaningless. In what follows, we shall develop a method for personalised pricing of the bond based on a feedback signal. The feedback signal will be designed so that aggregate levels of compliance satisfy some  constraint. Before proceeding we present two examples of social contracts and show how they lead to different policy choices.
The first example is akin to the traffic signal situation mentioned in the previous section. In this kind of application, a desirable policy might be the following: if the agent behaves, their token is returned, otherwise they lose their token. The second example concerns an agent who enters and moves within a public building (such as an airport or a train station) where the social contract might represent a rule such as {\em keep your mask on}. Clearly, in this type of contract, if the agent does not remove their mask then all tokens are returned. But what should happen if the agent does remove their mask? One policy might be to issue and redeem tokens at discrete intervals of time. Participants who break the contract multiple times would pay the bond repeatedly; this would also incentivise those individuals who remove their mask to wear it again so as to avoid further penalty. An alternative policy would be for an agent to lose their tokens, but for the pricing algorithm to operate at discrete intervals without further loss of tokens. In this case the agent is incentivised to wear the mask again so as to keep their personalised price as low as possible. Below we summarise some policies that are of interest to us.\newline  

\begin{itemize}
\item {\bf Fixed penalty policy: } Before participating in the social scheme each agent deposits a certain amount of tokens, the amount being set by the controller. When the action is completed or when the agent exits the scheme, all tokens are returned in the event that they complied with rule $E$; otherwise no tokens are returned to the agent. In the latter case the pricing algorithm continues to adjust the price based on both the agents' level of compliance and that of the network.\newline

\item{\bf Adaptive penalty policies:} Initially each agent deposits a certain amount of tokens, the amount being set by the controller. The contract is reissued at every time-step. 
At each time step, compliant agents retrieve their tokens, and stake new ones to continue the activity. Non-compliant agents lose all their tokens every time they do not comply. At all time steps, the pricing algorithm continues to adjust the price based on both the agents' level of compliance and that of the network.\newline

\item{\bf Adaptive penalty policies with return:} Initially each agent deposits a certain amount of tokens, the amount being set by the controller. The contract is reissued at every time-step. 
At each time step, compliant agents retrieve their tokens, and stake new ones to continue the activity. Non-compliant agents lose all their tokens every time they do not comply. If an agent that previously lost a token starts complying again, they will retrieve a portion of the lost tokens. At all time steps, the pricing algorithm continues to adjust the price based on both the agents' level of compliance and that of the network.\newline

\item {\bf Event driven policies: } Initially each agent deposits a certain amount of tokens. Whenever the agent fails to comply with rule $E$ the tokens are lost; in order to keep participating in the scheme the agent needs to deposit more tokens. In this version of the scheme the amount of tokens that are required varies as a bond changes value over time (again a smart contract can easily take care of the update process).\newline
\end{itemize}

Clearly, these are just four possible policies that might be adopted by the issuer of a social contract, and many others are possible. Our main contribution in this paper is to develop a modeling and feedback control strategy to describe and enable a wide class of policies that include the four aforementioned ones.

\section{Mathematical Framework}
{As previously explained, we are interested in designing a feedback mechanism to avoid scenarios in which the value of the bond is either too low (leading to non-compliance) or too high (meaning that agents would not engage in the scheme for fear of losing their tokens). The issue of finding this value is the subject of this section. We will use typical elements of control theory in a stochastic environment where a large number of agents interact with one another and are subject to rule $E$.\newline

Accordingly  we consider $n$ agents and, for each of them,  we define dependent binary random variables $\{ M_i(k) \in \{0,1\} \}_{i=1}^n$, for discrete values of $k$, such that

\begin{equation}
  \P(\text{$i$ complies with rule $E$ at time $k$}) = \P(M_i(k) = 1)
 \label{eq: prob}
\end{equation}

Moreover, we assume that the probability of these events is entirely dependent on a constant $q_i$, which represents the proclivity of each agent to comply with rules, and two control variables, $C(k)$, $c_i(k)$. The variable $C(k)+c_i(k)$ represents the value of the token bond staked by agent $i$ at time-step $k$. The combination $q_i + C(k)+c_i(k)$ determines the likelihood that agent $i$ will comply with the rule at time-step $k+1$. Then, (\ref{eq: prob}) can be expressed as

\begin{equation}
    \P(M_i(k+1) = 1) = p\left(q_i + C(k)+c_i(k)\right)
\end{equation}

with $p: \mathbb{R} \longrightarrow [0,1]$ being a monotone increasing function (which is used to bind the probability between 0 and 1). $C(k)$ and $c_i(k)$ represent, respectively, a global and an individual feedback signal whose purpose is to regulate the behaviour of each agent so as to achieve the desired level of compliance. However, due to the fact that the agents use a Distributed Ledger as a medium of communication, they can only access past levels of compliance $\{ M_i(k - m) \}_{i=1}^n$, where $m$ is a delay in the measurements caused by POW, synchronization across ledgers and verification time \cite{IoTPaper1}. 
Accordingly, we consider the following control laws, $\forall k \in \mathbb{N}$ and $\forall
 i \in \{1,\dots, n\}$,

\be\label{recur1}
C(k+1) &=& C(k) + \alpha \, \left(Q^* - n^{-1} \sum_{i=1}^n M_i(k - m)\right)\nonumber \\
c_i(k+1) &=& c_i(k) + \beta \, \left(Q^* - \overline{M_i(k - m)}\right)
\ee
with $\alpha > 0$ and $\beta > 0$ being two constants, $Q^* \in [0,1]$ being the desired level of compliance and $\overline{M_i(k)}$ representing a windowed time average of the compliance of agent $i$, defined as

\be\label{def:M-bar}
\overline{M_i(k)} = (1 - \gamma) \, \sum_{j=1}^k \gamma^{k-j} \, M_i(j).
\ee
In this last expression,
the factor $(1 - \gamma)^{-1}$ plays the role of the length of the window for the average, with  $\gamma < 1$. 

Notice that the proposed framework is very flexible and it would be possible to employ more sophisticated control laws. In this paper, however, we limit ourselves to the study of a proportional action and the extension to more complex feedback loops will be the subject of a future work. \newline
The reason to use both a global and an individual control signal, as opposed to just an individual or a global one,  is that these two feedback signals achieve different complementary goals:\newline

\begin{itemize}

\item \emph{Fairness:}  Due to differences in individual behaviour, some agents are going to comply with rules less than others. This means that if only a global shared signal were used to control the behaviour of multiple agents, the signal would be driven up by the behaviour of the least complying users, and this would result in an unfair price for the most \emph{virtuous} agents. On the other hand, the introduction of a personalised cost ensures that individuals are going to be priced according to their own behaviour (e.g., the less you comply the more you are going to pay, and vice versa).\newline

\item \emph{Distributed trading of compliance levels:} While the presence of a global cost is not necessary to guarantee the desired level of average compliance in the presented framework (because if every agent's compliance signal were equal to the target value $Q^*$ then the overall average compliance would be $Q^*$), its absence would make the system vulnerable to the repeated misbehaviour of malicious agents who  purposely attempt to drive down the compliance level. The introduction of a global signal ensures that the system is able to achieve the desired level of compliance even in the presence of this kind of disturbance. In effect, the global cost allows compliant agents to compensate for non-compliant ones. This aspect is further explored in Section \ref{sec: Simulations}.\newline

\item {\em Pricing attacks:} In traditional pricing models, even one nefarious agent could, in principle, drive up the cost for all agents simply by misbehaving.  In view of the previous comment, a natural concern is that similar effects might be possible in our schemes. Fortunately, in our scheme such attacks are not possible. Even though non-compliant agents may drive up the control signal $C(k)$ and hence drive up the cost of the bond, compliant agents will always receive their full deposit after expiration of a contract, leaving them unaffected by the increased price. On the other hand non-compliant agents would continue to lose tokens while they drive $C(k)$ to a high value. Furthermore, in the event that non-compliant agents prevent the desired levels of compliance being reached, $c_i(k)$ may in fact tend to zero for compliant agents, further rewarding their behaviour. \newline 

\end{itemize}

{\em Remark:} Before moving on, we want to point out that, since we are considering a DLT as the communication backbone of the whole architecture, the loss of a token is recorded on each agent's copy of the ledger. This means that everyone is aware at all times (subject to delay) of the level of compliance of every other actor in the scheme  (as the loss of the token by agent $i$ at time $k$ is recorded by the non-compliance value $M_i(k)=0$). Therefore, equations (\ref{recur1})  can be validated by each user individually. Also, notice that equations (\ref{recur1}) are consistent with the scenarios described in Section \ref{sec: Policy}.

\section{Theoretical Analysis}

In this section we provide a theoretical analysis of the convergence properties of the stochastic processes $\{C(k), c_i(k), \overline{M_i(k)}\}$.
Recall that the cost functions evolve over one time step as
\be
C(k+1) &=& C(k) + \alpha \,\left(Q^* - \frac{1}{n} \, \sum_{i=1}^n M_i(k - m)\right) \nonumber \\
c_i(k+1) &=& c_i(k) + \beta \, \left(Q^* - \overline{M_i(k - m).}\right)  \label{recur12}
\ee
and the time averaged variable $\overline{M_i(k)}$ satisfies the recursion formula
\be\label{recur1234}
\overline{M_i(k+1)} = &&\hskip-0.25in  \gamma \, \overline{M_i(k)} + (1 - \gamma) \, M_i(k+1)
\ee
We will use the following ReLU-type choice for the agent probability function:
\bee
p(x) =  \mid \{0, 1, x\} = \begin{cases}
    x & \text{if } x \in  (0,1)\\
    0 & \text{if } x \leq 0\\
    1 & \text{if } x \geq 1
\end{cases}.
\eee
Our main result is a convergence of probability for the time averaged compliance variables
$\overline{M_i(k)}$ around the target value $Q^*$ in the regime
\bee
\alpha << \begin{cases} \hskip0.1in \beta \cr 1 - \gamma \end{cases} \hskip-0.18in\Bigg\} << 1
\eee
Accordingly we introduce a small parameter $\epsilon$ and a window size parameter $w$,
and we allow $\alpha, \beta, \gamma$ to scale with $\epsilon, w$ as follows:
\be\label{def:scaling}
\alpha = \epsilon^{3/2} \, \alpha_0, \quad
\beta = \epsilon \, \beta_0 \, w, \quad
1 - \gamma = \epsilon \, w
\ee
where $\alpha_0, \beta_0$ are fixed constants.
Our results will hold for $\epsilon$ sufficiently small.
We will also assume that the parameter $\beta_0$ satisfies the following condition:
\be\label{beta0}
0 < \beta_0 \le 1
\ee
Since the recursion equations (\ref{recur12}) involve the time delay $m$, it is necessary to include initial conditions
\be\label{C-init}
\mathcal{C} =\{c_i(k-m),C(k-m),M_i(k-m),\overline{M_i(k-m)}\}
\ee
for all $k=0,\dots,m$ and $i=1,\dots,n$. We say $\mathcal{C}$ is normalized if
for all $k=0,\dots,m$ and $i=1,\dots,n$,
\be\label{init.norm1}
q_i +C(k-m) + c_i(k-m)  \in  [0,1], \nonumber \\
 M_i(k-m) \in \{0,1\}, \, \overline{M_i(k-m)} \in  [0,1].
\ee

\medskip
\begin{thm}\label{thm1}
For any $Q^* \in (0,1)$,
suppose that the processes $\{C(k), c_i(k), \overline{M_i(k)}\}$ satisfy the recursion equations (\ref{recur12}) and (\ref{recur1234})
and the parameters $\alpha, \beta, \gamma$ satisfy the relations (\ref{def:scaling}) and (\ref{beta0}). 
Then there are positive constants $B, \epsilon_0$ such that for all $\epsilon \le \epsilon_0$,
all normalized initial conditions $\mathcal{C}$ satisfying (\ref{init.norm1}), and all $\delta > 0$ and $i \in \{1,\dots,n\}$,
\be\label{thm1:eq1}
\limsup_{k \rightarrow \infty}  \P\left(|\overline{M_i(k)} - Q^*| > \delta\right)  \le  B \, \delta^{-2} \, \epsilon.
\ee
Taking $\delta = \epsilon^{1/3}$ we obtain that for all $\epsilon \le \epsilon_0$,
\be\label{thm1:eq2}
\limsup_{k \rightarrow \infty} \P\left(|\overline{M_i(k)} - Q^*| > \epsilon^{1/3}\right)  \le  B \, \epsilon^{1/3}.
\ee
\end{thm}

\medskip
\par\noindent{\em Remark }: 
The result of Theorem \ref{thm1} says that with high probability, the time average compliance level $\overline{M_i(k)}$ 
will be found close to the target value $Q^*$ for every agent $i$, for all $k$ sufficiently large. 
Note that the recursion equations (\ref{recur12}) and (\ref{recur1234}) can be recast as stochastic approximation equations
with constant stepsize $\epsilon$. As such, our model does not fit into the Robbins-Monro setting \cite{Rob-Mon} where the stepsize decreases to zero
as $k \rightarrow \infty$, and so we cannot expect to derive concentration-type results showing that the sequence $\overline{M_i(k)}$ 
converges to a neighborhood of  $Q^*$ as $k \rightarrow \infty$. Instead Theorem \ref{thm1}
is a statement about the long-term distribution of $\overline{M_i(k)}$, and shows that the
distribution is concentrated close to $Q^*$ for all $k$ sufficiently large.
Note that our proof also shows that with high probability the cost function $q_i+ C(k) + c_i(k)$
will be found close to the target value $Q^*$ for every agent $i$, for large $k$.

\medskip
\par\noindent{\em Remark }: 
Although we chose a specific form for the agent probability function $p$, our results can be extended to allow personalized
probability functions $\{p_i\}$ for each agent,  where in each case
$p_i : \mathbb{R} \rightarrow [0,1]$ is a non-decreasing uniformly Lipschitz function.

\subsection{Proof of Theorem \ref{thm1}}
In order to prove Theorem 1 we will first construct deterministic difference equations for each agent which approximate
the stochastic recursion equations (\ref{recur12}) and (\ref{recur1234}). We will then show that the solutions of the
stochastic and deterministic equations remain close in $L^2$ norm as $k \rightarrow \infty$, and that the solution of the
deterministic difference equations converges to a small neighborhood of the target value $Q^*$.
The $L^2$ norm will be written as
\be\label{def-L2}
\| X \|_{L^2} = \left(\E[\|X\|^2]\right)^{1/2}
\ee

\medskip

Let ${\cal F}(k)$ be the $\sigma$-algebra generated by the random variables
$\{M_i(j) : 1 \le i \le n, \, 1 \le j \le k\}$. Then conditioned on ${\cal F}(k)$, the Bernoulli random variables 
$\{M_i(k+1)\}$ ($i=1,\dots,n$) are independent with the distribution
\be\label{def:M}
\P(M_i(k+1)=1 \,|\, {\cal F}(k)) = p(q_i+ C(k) + c_i(k)).
\ee
We also define the variables $\{\xi_i(k)\}$ by
\bee
\xi_i(k) &=& M_i(k) - p(q_i+ C(k-1) + c_i(k-1))
\eee
It follows that 
\bee
\E[\xi_i(k+1) \,|\, {\cal F}(k)] = 0
\eee
and therefore $\{\xi_i(k)\}$ form a martingale difference sequence with respect to ${\cal F}(k)$.
For future use we also note  the bound
\be\label{xi-bd}
\E\left[\xi_i(k+1)^2 \,|\, {\cal F}(k)\right] = {\rm VAR}[M_i(k+1) \,|\, {\cal F}(k)] \le \frac{1}{4}
\ee
The recursion formula (\ref{recur1234}) can be written as
\be
\overline{M_i(k+1)} = &&\hskip-0.25in  \gamma \, \overline{M_i(k)} + (1 - \gamma) \, M_i(k+1) \nonumber \\
= &&\hskip-0.25in \overline{M_i(k)} + (1 - \gamma) \, \xi_i(k+1) \nonumber \\
&& \hskip-0.4in + (1 - \gamma) \,  \left[p(q_i+ C(k) + c_i(k)) - \overline{M_i(k)}\right] \label{recur123}
\ee

\medskip
We now choose some index $i$ corresponding to a particular agent: this index will be fixed throughout the proof.
We define for each $k$
\be\label{def-Y}
Y_1(k) &=& \overline{M_i(k)} \\
Y_2(k) &=& q_i + C(k) + c_i(k)
\ee
and write $Y(k) = (Y_1(k), Y_2(k))^T \in \mathbb{R}^{2}$. We define the function 
$h =(h_1,h_2): \mathbb{R}^{2} \rightarrow \mathbb{R}^{2}$ as follows:
\be
h_1(y) &=& w \, \left(p(y_{2}) - y_1\right), \\
h_{2}(y) &=&  \beta_0 \, w \, \left(Q^* - y_1 \right), \\
\ee
Using these definitions we can rewrite the recursion equation (\ref{recur123}) as follows:
\be
Y_1(k+1) &=& Y_1(k) + (1 - \gamma) \, \xi_i(k+1) \nonumber \\
&&  + (1 - \gamma) \,  \left[p(q_i+ C(k) + c_i(k)) - \overline{M_i(k)}\right] \nonumber \\
&=& Y_1(k) + \epsilon \, w \, \xi_i(k+1) + \epsilon \,  h_1(Y(k)) \label{recurY1}
\ee
and similarly from (\ref{recur12}):
\be
Y_2(k+1) = Y_2(k) + \epsilon \,  h_2(Y(k-m)) + \epsilon^{3/2} \, G(k) \label{recurY2}
\ee
where
\be
G(k) = \alpha_0 \, \left(Q^* - \frac{1}{n} \, \sum_{i=1}^n M_i(k - m)\right) \label{def:G}
\ee
The initial conditions $\mathcal{C}$ as in  (\ref{C-init})
define initial conditions $\{Y(k-m), \, k=0,\dots,m\}$ for (\ref{recurY1}) and (\ref{recurY2}).

\medskip
We next introduce deterministic difference equations by dropping terms from (\ref{recurY1}) and (\ref{recurY2}):
\be
y_1(k+1) &=& y_1(k) + \epsilon \,  h_1(y(k)) \nonumber \\
y_2(k+1) &=& y_2(k) + \epsilon \,  h_2(y(k-m)) \label{recury2}
\ee
where $y(k) = (y_1(k),y_2(k)) \in \mathbb{R}^2$. 
It is evident that $(Q^*,Q^*)$ is the unique fixed point of this system.
The initial conditions for (\ref{recury2}) are $\{y(k-m), \, k=0,\dots,m\}$.
Elementary estimates (using $0 \le p(y_2) \le 1$) show that if $y_1(0) \in [0,1]$ then $y_1(k) \in [0,1]$ for all $k$,
and therefore $|h_1(y(j))| \le w$ for all $j \ge 0$. Assuming also that $y_1(k-m) \in [0,1]$ for all $k=0,\dots,m$
it follows that $|h_2(y(j))| \le \beta_0 w \le w$ for all $j \ge -m$, and therefore
\be\label{bd-ydiff}
\| y(j+1) - y(j) \| \le 2 \,w \, \epsilon \quad \text{for all $j \ge 0$}.
\ee
We will say that the initial condition $\mathcal{Y} = \{y(k-m)\}$ is {\em bounded at level $M$} if 
\be
&& \hskip-0.25in y_1(k-m) \in [0,1], \,\, |y_2(k-m)| \le M \,\, \forall \,\, k=0,\dots,m. \nonumber \\
&& \label{norm-def2}
\ee
Note that (\ref{bd-ydiff}) holds for any initial condition $\mathcal{Y}$ bounded at level $M$, for any $M \ge 0$.

\medskip
The proof of Theorem \ref{thm1} proceeds by first bounding the difference between the stochastic sequence $Y(k)$
and the deterministic sequence $y(k)$. 
We assume that  $Y(k)$ and $y(k)$ share the same initial conditions $\mathcal{Y}$.
Then 
\be
Y_1(k) - Y_1(0) &=& \epsilon \, \sum_{j=0}^{k-1} h_1(Y(j)) + \epsilon \, w \, \sum_{j=0}^{k-1} \xi_i(j+1) \nonumber \\
y_1(k) - y_1(0) &=& \epsilon \, \sum_{j=0}^{k-1}   h_1(y(j))
\ee
Since $Y_1(0)=y_1(0)$ we get
\be
Y_1(k) - y_1(k)& =& \epsilon \, \sum_{j=1}^{k-1} h_1(Y(j)) - h_1(y(j)) \nonumber \\
&& \hskip0.2in + \epsilon \,w\, \sum_{j=0}^{k-1} \xi_i(j+1)  \nonumber 
\ee
The function $h_1$ is uniformly Lipschitz on $\mathbb{R}^{2}$, and
\bee
| h_1(y) - h_1(z) |  \le  2 \, w \, \| y - z \| \quad \text{for all $y,z \in \mathbb{R}^{2}$}.
\eee
Therefore
\be\label{Y1bd}
&& \hskip-0.6in \| Y_1(k) - y_1(k) \|_{L^2} \nonumber \\
&& \le 2 \, \epsilon \, w \, \sum_{j=1}^{k-1} \| Y(j) - y(j) \|_{L^2} + \epsilon \, w \, \frac{1}{2} \, k^{1/2} \nonumber \\
&&
\ee
where we used martingale orthogonality and (\ref{xi-bd}) to deduce that
\bee
&&\hskip-0.4in \|\sum_{j=0}^{k-1} \xi_i(j+1)\|^{2}_{L^2} = \sum_{j=0}^{k-1} \E\left[\xi_i(j+1)^2\right]  \le  \frac{k}{4}
\eee
Similarly
\bee
Y_2(k) - y_2(k)& =& \epsilon \, \sum_{j=0}^{k-1} h_2(Y(j-m)) - h_2(y(j-m)) \\
&& \hskip0.2in + \epsilon^{3/2} \, \sum_{j=0}^{k-1} G(j)
\eee
Noting that $|G(k)| \le \alpha_0$ and $|h_2(y)-h_2(z)| \le w \|y-z\|$
we deduce
\be
&& \hskip-0.3in \| Y_2(k) - y_2(k) \|_{L^2} \hskip-0.1in \nonumber \\
&&  \hskip-0.5in \le   \epsilon \, w \, \sum_{j=m+1}^{k-1} \| Y(j-m) - y(j-m) \|_{L^2} 
+ \epsilon^{3/2} \, \alpha_0 \, k \label{Y2bd}
\ee
(where we used $Y(j-m) = y(j-m)$ for all $j=0,\dots,m$).
Combining (\ref{Y1bd}) and (\ref{Y2bd}) gives
\be
\| Y(k) - y(k) \|_{L^2} &\le & 3\, \epsilon  \, w \, \sum_{j=1}^{k-1} \| Y(j) - y(j) \|_{L^2} \nonumber \\
&& + \epsilon \, w \, \frac{1}{2} \, k^{1/2} + \epsilon^{3/2} \, \alpha_0 \, k
\ee
Therefore by applying the discrete Gr\"onwall inequality we deduce that for all $k$
\be\label{Ybd3}
\sup_{0 \le j \le k} \| Y(j) - y(j) \|_2 \le \epsilon^{1/2} \, \left(\frac{(\epsilon k)^{1/2} \, w}{2} + \epsilon \, k \, \alpha_0\right) \, e^{3 \, w \epsilon k}
\ee

\medskip
The bound (\ref{Ybd3}) implies that $Y(k)$ remains close to $y(k)$ over time intervals
$k \sim \epsilon^{-1}$, for any initial conditions shared by $Y$ and $y$. 
We next show that $y(k)$ converges to a neighborhood of $(Q^*,Q^*)$, using
some ideas and techniques from \cite{Borkar}, Chapter 9.

\begin{lemma}\label{lem2}
\noindent{\bf a)} 
For any $M \ge 0$, there are positive constants $B_1, \tau$ such that for all $\epsilon$ sufficiently small and
all initial conditions $\mathcal{Y}$ bounded at level $M$, the sequence $y(k)$ defined by
(\ref{recury2})  satisfies
\be\label{lem2:eq2}
\| y(\lfloor\tau \epsilon^{-1}\rfloor) - (Q^*,Q^*) \| \le \frac{1}{2} \, \| y(0) - (Q^*,Q^*) \| + B_1 \, \epsilon
\ee
\noindent{\bf b)} There is $M_2 > 2$ such that for all $\epsilon$ sufficiently small and 
all initial conditions $\mathcal{Y}$ bounded at level $1$, the sequence $Y(k)$ defined by
(\ref{recurY1}), (\ref{recurY2}) satisfies
\be
0 \le Y_1(k) \le 1, \,\, - M_2 \le Y_2(k) \le  M_2 \quad \text{for all $k \ge 0$} \label{lem2:eq1b}
\ee
\end{lemma}

Lemma \ref{lem2} will be proved in the Appendix. 
We can now complete the proof of Theorem \ref{thm1}.
As a convenient shorthand we will write
\be
Y(k:k+l) = \{Y(k),Y(k+1),\dots,Y(k+l)\}
\ee
to denote a sequence of successive states of the process. Also, we write $\E_{\mathcal{Y}}$ to denote expected value
with initial condition $\mathcal{Y}$ for the system (\ref{recurY1}),(\ref{recurY2}).
We define $J = \lfloor\epsilon^{-1}\tau\rfloor$ where the number $\tau$ is the value defined in Lemma \ref{lem2}(a) for $M = M_2$.
Since the distribution of the increment of the system (\ref{recurY1}), (\ref{recurY2}) depends only on the previous $m+1$ states, 
it follows that for any $\mathcal{Y}$
and integer $k \ge 0$ we have
\be\label{cond-bd2}
&&\hskip-0.6in \E\left[\|Y(k+J) - (Q^*,Q^*)\|^2 \,|\, Y(k-m : k) = \mathcal{Y}\right] \nonumber \\
&& \hskip0.2in = \E_{\mathcal{Y}}\left[\|Y(J) - (Q^*,Q^*)\|^2\right]
\ee
We define $y(l)$ as the solution of (\ref{recury2}) with initial conditions $\mathcal{Y}$,
and define
\be
B_2 =  \left(\frac{\tau^{1/2} \, w}{2} +\tau \, \alpha_0\right) \, e^{3 w \tau}
\ee
The bound (\ref{Ybd3}) implies that
\be\label{bd7}
\E_{\mathcal{Y}}\left[\|Y(J) - y(J)\|^2\right] \le \epsilon \, B_2^2
\ee
Furthermore from Lemma \ref{lem2}(b) we may assume that $\mathcal{Y}$ is normalized at level $M_2$,
and so from Lemma \ref{lem2}(a) we deduce that
\be\label{bd8}
\|y(J) - (Q^*,Q^*)\|^2 \le \frac{\sqrt{2}}{4} \, \| y(0) - (Q^*,Q^*) \|^2 + B_1'^2 \, \epsilon^2
\ee
where $B_1'^2=(2 + \sqrt{2}) \, B_1^2$.
Combining (\ref{bd7}) and (\ref{bd8}) gives
\be
&& \hskip-0.2in \E_{\mathcal{Y}}\left[\|Y(J) - (Q^*,Q^*)\|^2\right] \nonumber \\
&& \hskip-0.3in \le \sqrt{2} \, \|y(J) - (Q^*,Q^*)\|^2 + (2 + \sqrt{2}) \, \E_{\mathcal{Y}}\left[\|Y(J) - y(J)\|^2\right]  \nonumber \\
&& \hskip-0.3in \le  \frac{1}{2} \, \| y(0) - (Q^*,Q^*) \|^2 + \epsilon \, B_3^2
\ee
where (using $\epsilon < 1$)
\be
B_3^2 = (2 + \sqrt{2}) \, B_2^2 + \sqrt{2} \, B_1'^2
\ee
Substituting into (\ref{cond-bd2}) gives
\be
&&\hskip-0.6in \E\left[\|Y(k+J) - (Q^*,Q^*)\|^2 \,|\, Y(k-m : k)\right] \nonumber \\
&& \hskip0.2in \le \epsilon \, B_3^2 + \frac{1}{2} \, \| Y(k) - (Q^*,Q^*) \|^2
\ee
and therefore
\be\label{recursive-bd1}
&&\hskip-0.6in \E\left[\|Y(k+J) - (Q^*,Q^*)\|^2\right] \nonumber \\
&& \hskip0.2in \le \epsilon \, B_3^2 + \frac{1}{2} \, \E\left[\| Y(k) - (Q^*,Q^*) \|^2\right]
\ee
For any $k \ge 1$ we write $k = n J + \wh{k}$ where $J = \lfloor\epsilon^{-1}\tau\rfloor$,
$n \ge 0$ and $0 \le \wh{k} \le J-1$.
Then applying (\ref{recursive-bd1}) recursively and using the bound (\ref{lem2:eq1b}) gives
\be\label{recursive-bd2}
&&\hskip-0.4in \E\left[\|Y(k) - (Q^*,Q^*)\|^2\right] \nonumber \\
&& \hskip0.2in \le 2 \, \epsilon \, B_3^2 + 2^{-n} \, \E\left[\| Y(\wh{k}) - (Q^*,Q^*) \|^2\right] \nonumber \\
&& \hskip0.2in \le 2 \, \epsilon \, B_3^2 + 2^{-n} \,  \left( 1 + (1 + M_2)^2\right).
\ee
Since $n \rightarrow \infty$ as $k \rightarrow \infty$ it follows that
\be\label{limsup-E}
\limsup_{k \rightarrow \infty} \E\left[\|Y(k) - (Q^*,Q^*)\|^2\right] \le 2 \, \epsilon \, B_3^2
\ee
Finally using Markov's inequality we get
\be\label{limsup-P}
\hskip-0.2in \P\left(|\overline{M_i(k)} - Q^*| > \delta\right)\hskip0.05in && \hskip-0.3in = P\left(|Y_1(k) - Q^*| > \delta\right) \nonumber \\
&& \hskip-0.3in \le \P\left(\| Y(k) - (Q^*,Q^*) \| > \delta\right) \nonumber  \\
&& \hskip-0.3in \le \delta^{-2} \, \E\left[\|Y(k) - (Q^*,Q^*)\|^2\right]
\ee
and (\ref{thm1:eq1}) then follows by combining (\ref{limsup-E}) and (\ref{limsup-P}), with
\be
B = 2 \, B_3^2\newline
\ee

\emph{Remark: }As a final remark, we comment on the relationship between the proposed control algorithms, and the class of 
algorithms that are commonly known as {\em consensus algorithms}. Consensus problems typically consider partial exchange of information between agents with asynchronous updates, and look for local update rules so that all agents converge to a common value (or variations thereof).  However, our problem is much simpler; agents have access to a global signal and full information from other agents via the distributed ledger. So even though the agents do agree asymptotically (if this is the pricing strategy), as in a consensus problem, there is no real explicit asynchronous exchange of partial information; rather all agents have full access to all agents' information. Thus, the most natural formulation is the one we chose; that of a conventional feedback regulation problem with the main theoretical contribution being that of convergence of the feedback loop {\em in probability}.

\section{Simulations}
\label{sec: Simulations}
In this section we provide simulations for the proposed control scheme, to show the effectiveness of our approach and to highlight the need for both a global and individual signal.\newline
Specifically, we consider the following scenarios:\newline
\begin{itemize}
\item[I] A scenario where only the global signal $C(k)$ is used to regulate the behaviour of each agent (i.e., $c_i(k) = 0, \, \forall k, \, \forall i$);
\item[II] A scenario where both global and individual signals are used to regulate the behaviour of each agent;
\item[III] A scenario where only the individual signal is used  to regulate the behaviour of each agent and a subset of individuals $D \subset \{1,\dots,n\}$ refuses to comply with rule $E$ (i.e., $\P (M_i(k) = 1) = 0, \forall k, \forall i \in D$);
\item[IV] A scenario where both the global and the individual signals are used  to regulate the behaviour of each agent and a subset of individuals $D \subset \{1,\dots,n\}$ refuses to comply with rule $E$.
\item[V] The same as scenario II but with increasing values of the delay $m$.
\item[VI] To analyse the robustness of the system in a more realistic scenario, where each agent is connected to a wireless network, we consider random delays (i.e., we allow the values of $m$ to be drawn according to a distribution) and the random event that at each time step an agent might experience disconnection or package drops. 
\end{itemize}

In all scenarios we set $n=1000$, $\alpha =0.025$, $\beta =0.1$, $Q^*=0.85$, $\gamma = 0.95$, $C(0) = 0$, $c_i(0) = 0, \, \forall i \in \{1,\dots, n\}$ and for scenarios I-IV, we consider $m = 3$. Base compliance levels $q_i$ are sampled from the uniform distribution in $[0.1, 0.35]$. Moreover, due to the stochastic nature of the system, for each scenario we perform 150 Monte Carlo simulations and we average over the obtained realizations in order to obtain statistically meaningful results. Simulations are performed on Matlab 2020a and Python 3.9.9.\newline

\subsection{Scenarios I and II}
Figures \ref{fig: global} and \ref{fig: individual} show the results of the simulations in scenario I and II. Even by visual inspection it is clear that the lack of an individual signal to control the agents' behaviour leads to unfair results: while the overall compliance converges to $Q^*$, this is achieved at the expense of the users that would behave better under normal circumstances (i.e., the users with larger $q_i$), that are forced to comply with a higher probability than $Q^*$ in order to compensate for the behaviour of the less compliant agents. Of course, this is undesirable and the use of the personalised cost, as shown in Figure \ref{fig: individual}, tackles this problem by adjusting the individual price depending on the past behaviour of each agent.
\subsection{Scenarios III and IV}
While in Scenarios I and II we explored how the lack of an individual cost leads to unfair results, it is less clear why a global cost is needed at all. In fact, equations (\ref{recur1}) show that, when $C(k)$ is set to zero, the personalised control signals would be sufficient to drive the average behaviour to the desired level of compliance. Nevertheless, without the global cost, the system might fail to achieve the desired target for compliance in scenarios when for some reason a certain number of agents fail to comply repeatedly with rule $E$. This could be due to malfunctions or malicious behaviour. This is highlighted in Figures \ref{fig: indOnly10}, related to scenario III, where 10 \% of the agents, for $k \leq 100$ does not comply with rule $E$ and the system is not able to achieve the desired level of compliance $Q^*$. On the other hand, in scenario IV, shown in Figures \ref{fig: Glob10}, it is possible to see that the presence of the global signal corrects this disturbance, thus making the system more robust to malfunctions and malicious behaviour (of course, the drawback is that the honest agents will have to comply more in order to compensate for the misbehaviour of the non compliant users).
\subsection{Scenario V}
As per the last set of simulations we show the behaviour of the system for increasing values of $m$. Figures \ref{fig: delay 10}-\ref{fig: delay 50}  show that while for values $m = 10$ and $m = 15$, the system maintains stability and the distributions of the average compliance and of each individual compliance accumulates around $Q^*$, for larger values , such as $m = 25$ and $m = 50$ the system ends up oscillating without ever reaching an equilibrium. This shows that the delay introduced by the DLT represents a crucial design parameter and that, for fast-paced applications, the choice of an architecture that allows quick approvals is of paramount importance. Of course, the results of this last scenario, are not meant to represent how a group of agent would behave, in the aforementioned circumstances. This scenario is merely showing a situation in which the control signal ends up failing due to the presence of the delay. 
{\color{black}\subsection{Scenario VI}
In our final set of experiments we are interested in the behaviour of the system when agents are connected to a network subject to heterogenous delays, such as a wireless or 5G network. In this scenario the system will experience random delays due to latency and potential disconnections or packet loss. While simulations based on the use of tools such as ns3 \cite{ns3} would be a more accurate representation of reality in such situations, as we are interested in the impact of delay on the feedback control algorithms, we restrict ourselves to Monte-Carlo based Python simulations to evaluate the performance of our algorithms. Integration of our work into an ns3 environment will be the subject of future work.\newline

Specifically, in this scenario we allow for $m$ to be drawn from a gaussian distribution with mean $\mu$ and standard deviation $\sigma$. This simulates the delay experienced by an agent connected to a wireless network (e.g., 3G or 4G). Moreover, at time step $k$, the level of compliance of agent $i$ will remain the same as its level of compliance at time step $k-1$, with probability $\eta$. This simulates the possibility that either the agent disconnects at time $k$ or the information never makes it to the controller.\newline

Accordingly, Figures \ref{fig: wifi5}-\ref{fig: wifi45} show the behaviour of the system for different values of $\mu, \sigma$ and $\eta$. More specifically, we set $\sigma = \mu/2$, $\eta = 0.025$ and we allow $\mu$ to vary between 5 and 45 with a step size of 10. Notice that the amount of time elapsed between time step $k$ and time step $k+1$ is application-dependent and therefore it would not be meaningful to provide $\mu$ expressed in physical time units.\newline 

As in scenario $V$, the system maintains stability and the distributions of the average compliance and of each individual compliance accumulates around $Q^*$, for small values of $\mu$, whereas for values of $\mu >50$ the system becomes unstable. The probability $\eta$ of an agent disconnecting or the information not reaching the controller does not affect the overall stability of the system. Notice that in scenario V, the system became unstable for values of $m>25$. Interestingly, drawing $m$ from a gaussian distribution, rather than being a constant value seems to increase the stability of the system. This aspect will be investigated in a future work.\newline

This shows that, similarly to the role of the delay introduced by a DLT, in scenario V, the delay introduced by wireless networks represents a crucial design parameter and that, for fast-paced applications, the choice of an architecture that allows quick approvals is of paramount importance. }

\begin{figure}

\includegraphics[width=1\columnwidth]{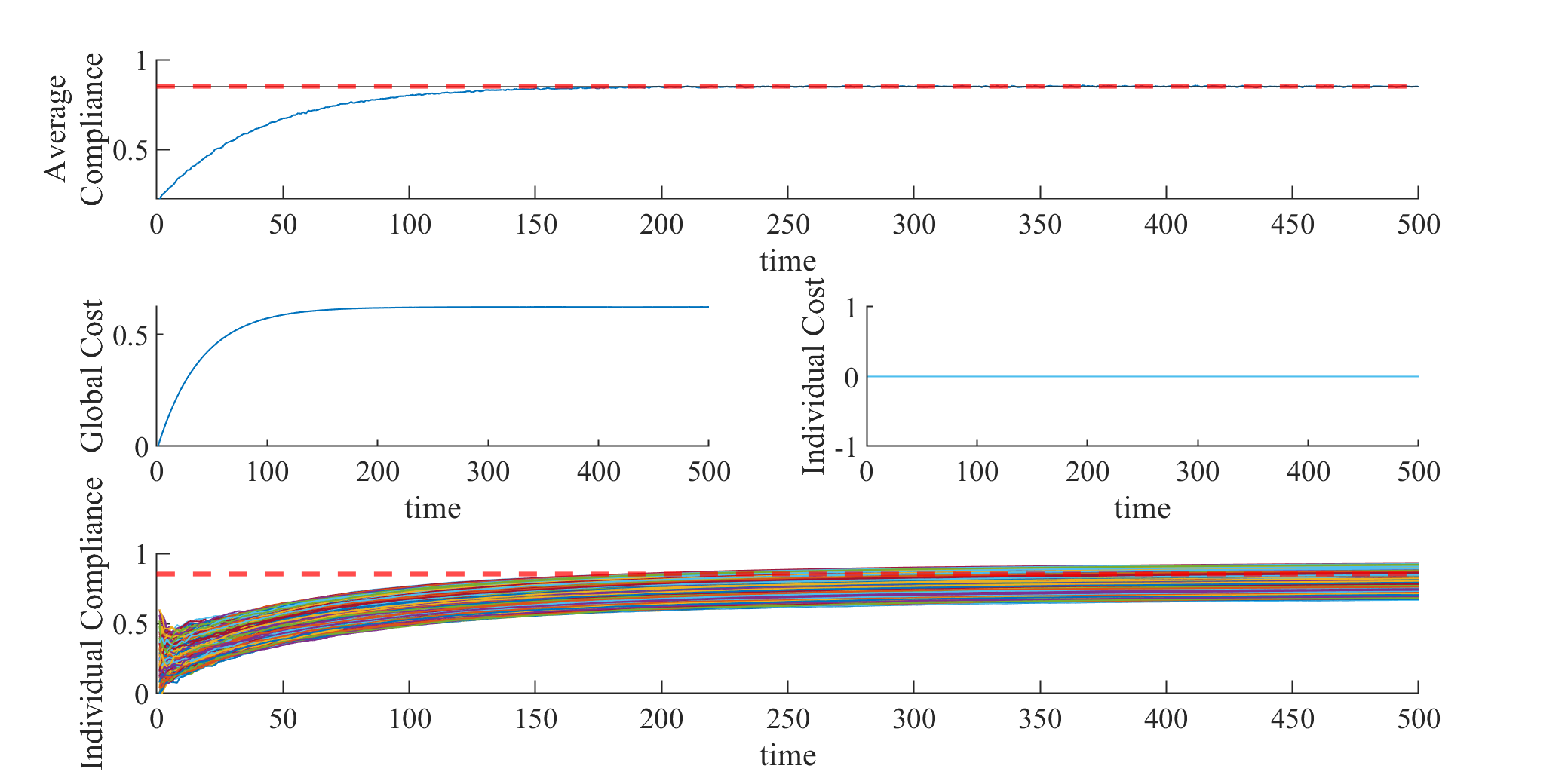}

\caption{Compliance control using only the global signal. Compliance is enforced but at the expense of some agents that have to comply more than others.}
\label{fig: global}
\end{figure}
\begin{figure}
\includegraphics[width=1\columnwidth]{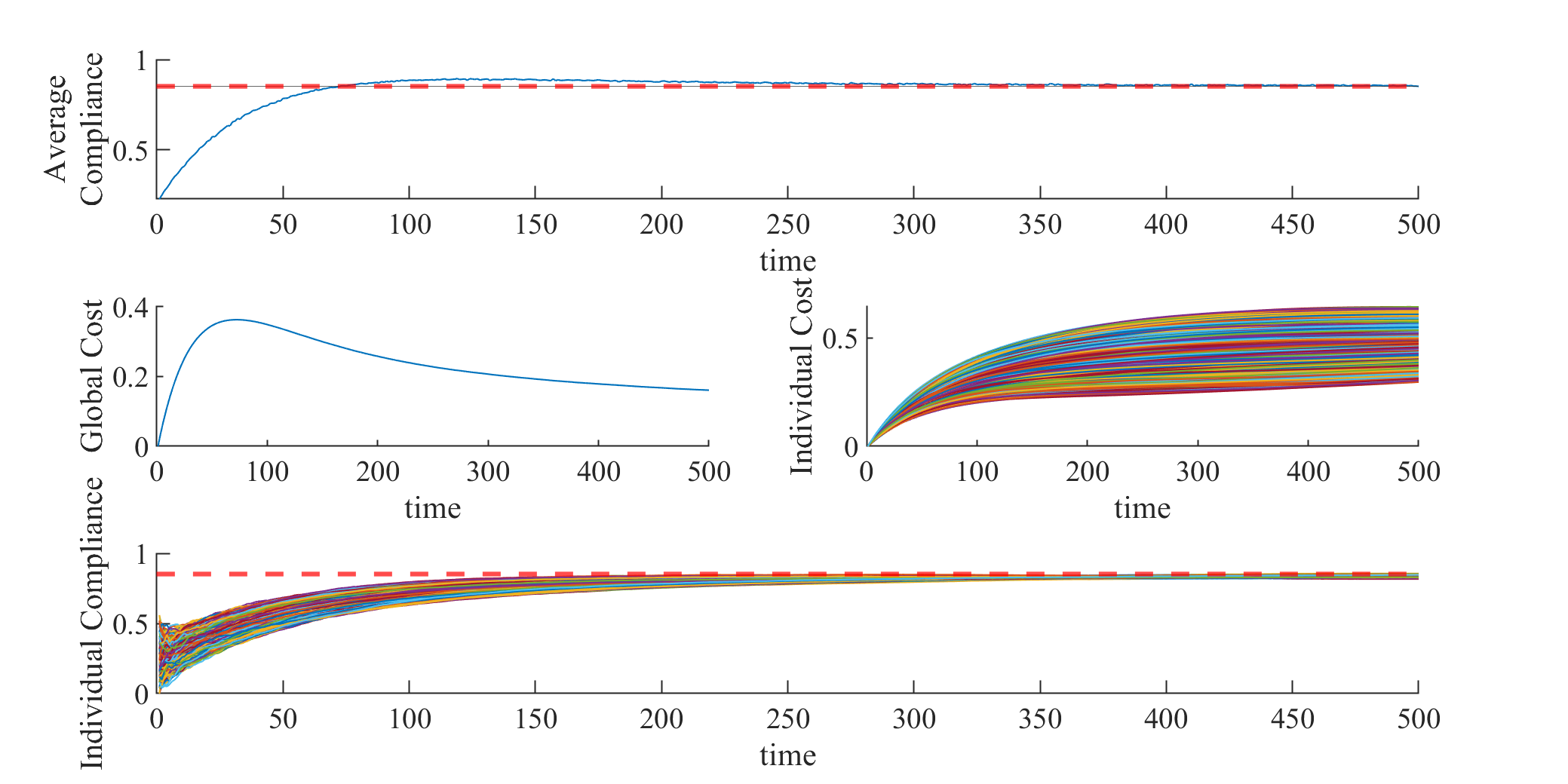}

\caption{Compliance control using both the global and the individual signals. Compliance is achieved and the use of personalised signals make it so that every agent contributes in a fair way.}
\label{fig: individual}
\end{figure}

\begin{figure}
\includegraphics[width=1\columnwidth]{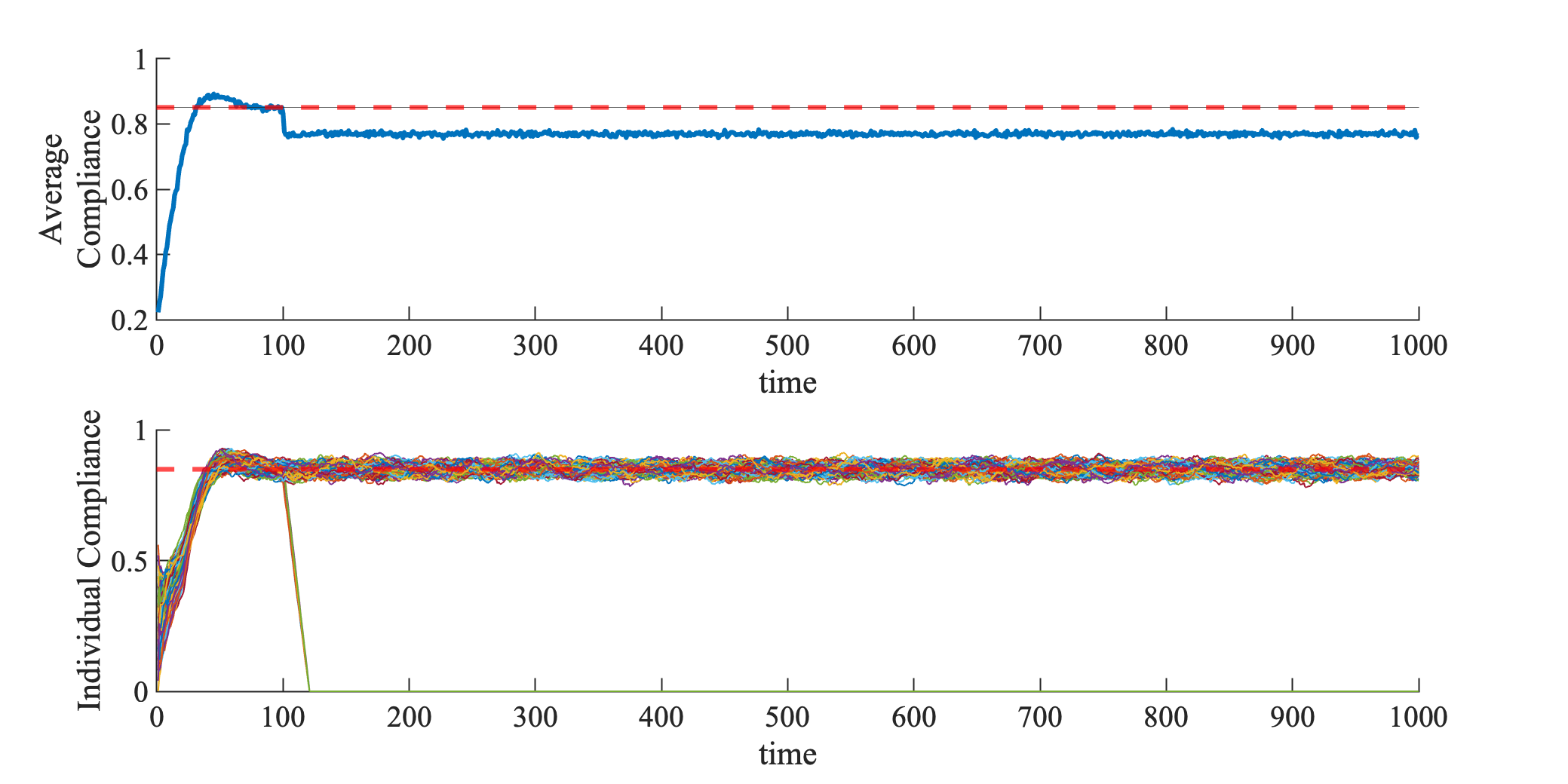}
\caption{Compliance control using only the individual signals. For $k \leq 100$ 10\% of the agents does not comply with rule $E$.}
\label{fig: indOnly10}
\end{figure}

\begin{figure}
\includegraphics[width=1\columnwidth]{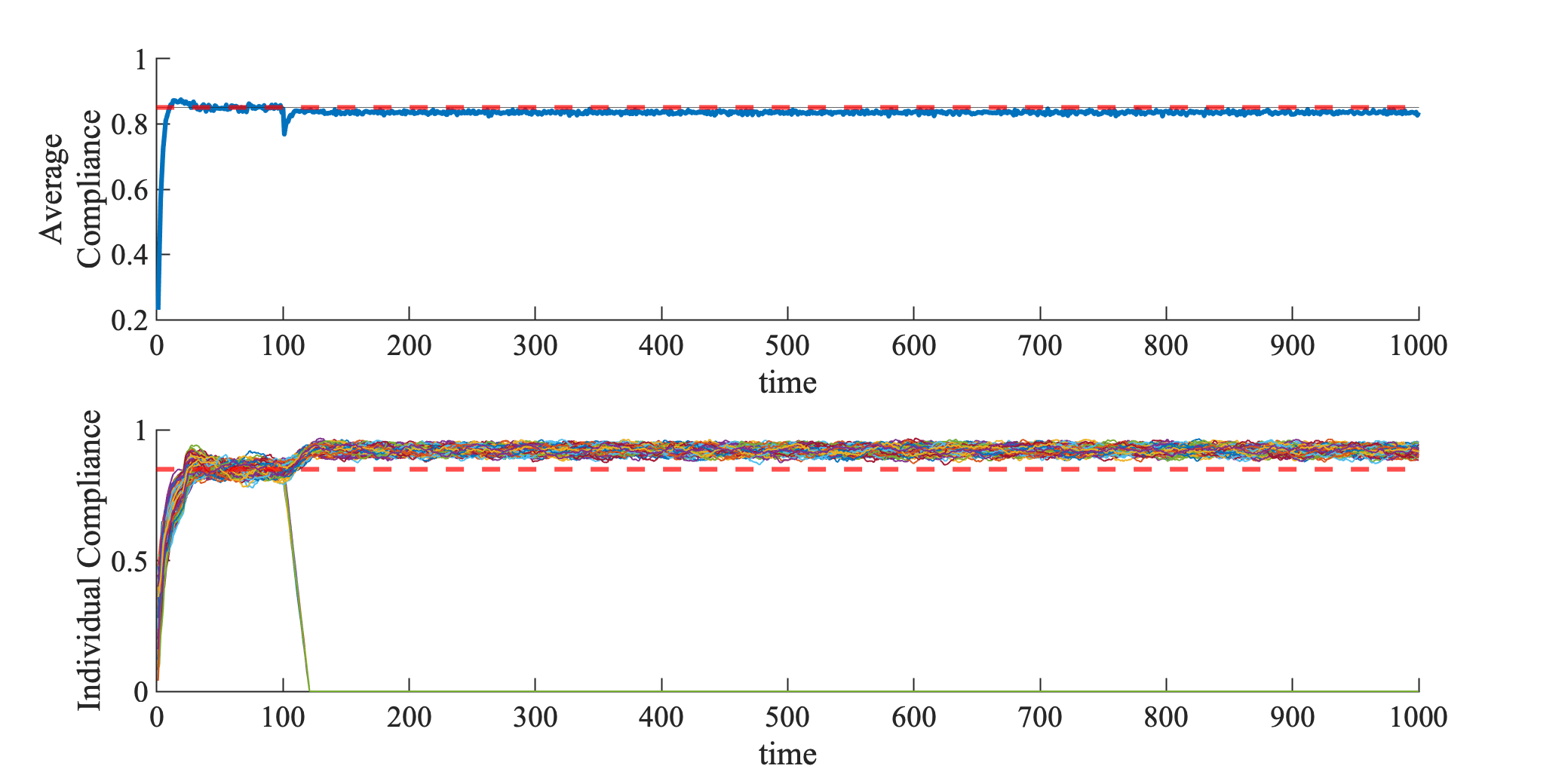}
\caption{Compliance control using both the global and the individual signals. For $k \leq 100$ 10\% of the agents does not comply with rule $E$.}
\label{fig: Glob10}
\end{figure}

\begin{figure}
\includegraphics[width=1\columnwidth]{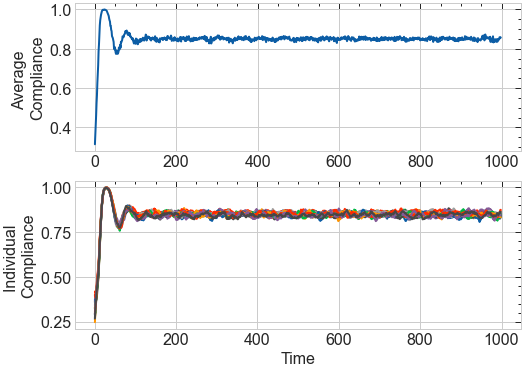}
\caption{Compliance control using both the global and the individual signals. Delay $m$ is set to 10. Compliance is achieved and the use of personalised signals make it so that every agent contributes in a fair way.}
\label{fig: delay 10}
\end{figure}

\begin{figure}
\includegraphics[width=1\columnwidth]{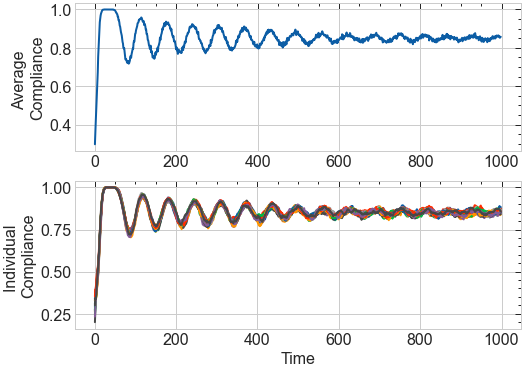}
\caption{Compliance control using both the global and the individual signals. Delay $m$ is set to 15. Compliance is achieved and the use of personalised signals make it so that every agent contributes in a fair way.}
\label{fig: delay 15}
\end{figure}

\begin{figure}
\includegraphics[width=1\columnwidth]{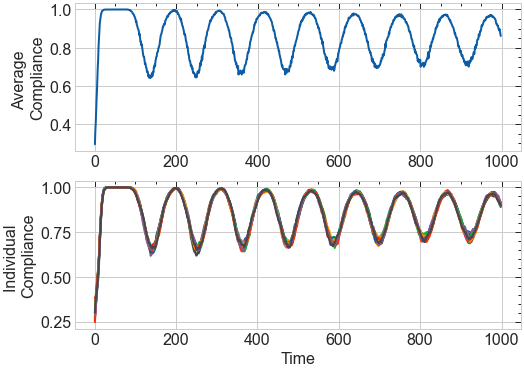}
\caption{Compliance control using both the global and the individual signals. Delay $m$ is set to 25. Compliance is not achieved and the system ends up oscillating.}
\label{fig: delay 25}
\end{figure}

\begin{figure}
\includegraphics[width=1\columnwidth]{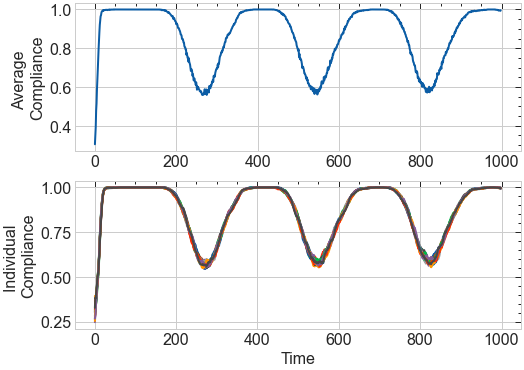}
\caption{Compliance control using both the global and the individual signals. Delay $m$ is set to 50. Compliance is not achieved and the system ends up oscillating.}
\label{fig: delay 50}
\end{figure}

\begin{figure}
\includegraphics[width=1\columnwidth]{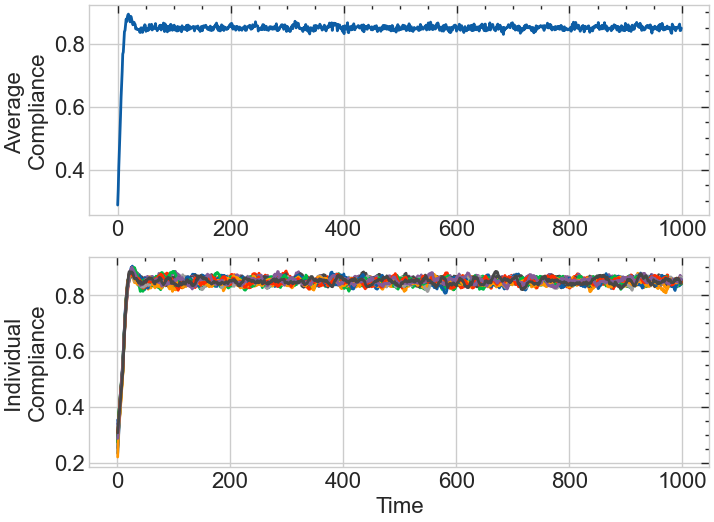}
\caption{Compliance control when the delay $m$ is a random variable sampled from a Gaussian distribution with parameters $\mu = 5$, $\sigma = 2.5$ and each agent has a chance of disconnecting from the network $\eta = 0.025$. Compliance is achieved and the use of personalised signals make it so that every agent contributes in a fair way.  }
\label{fig: wifi5}
\end{figure}

\begin{figure}
\includegraphics[width=1\columnwidth]{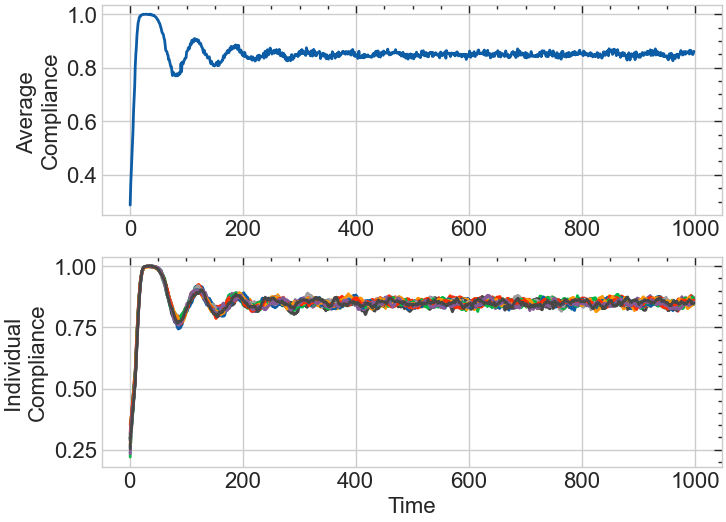}
\caption{Compliance control when the delay $m$ is a random variable sampled from a Gaussian distribution with parameters $\mu = 15$, $\sigma = 7.5$ and each agent has a chance of disconnecting from the network $\eta = 0.025$. Compliance is achieved and the use of personalised signals make it so that every agent contributes in a fair way.  }
\label{fig: wifi15}
\end{figure}

\begin{figure}
\includegraphics[width=1\columnwidth]{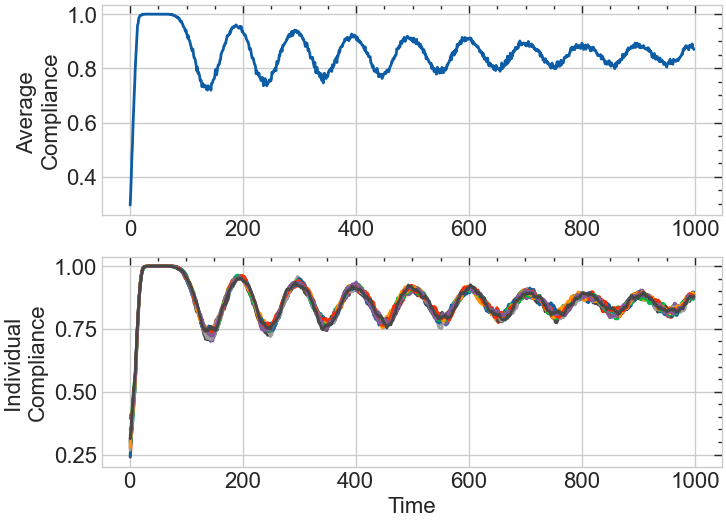}
\caption{Compliance control when the delay $m$ is a random variable sampled from a Gaussian distribution with parameters $\mu = 25$, $\sigma = 12.5$ and each agent has a chance of disconnecting from the network $\eta = 0.025$. Compliance is achieved and the use of personalised signals make it so that every agent contributes in a fair way.  }
\label{fig: wifi25}
\end{figure}

\begin{figure}
\includegraphics[width=1\columnwidth]{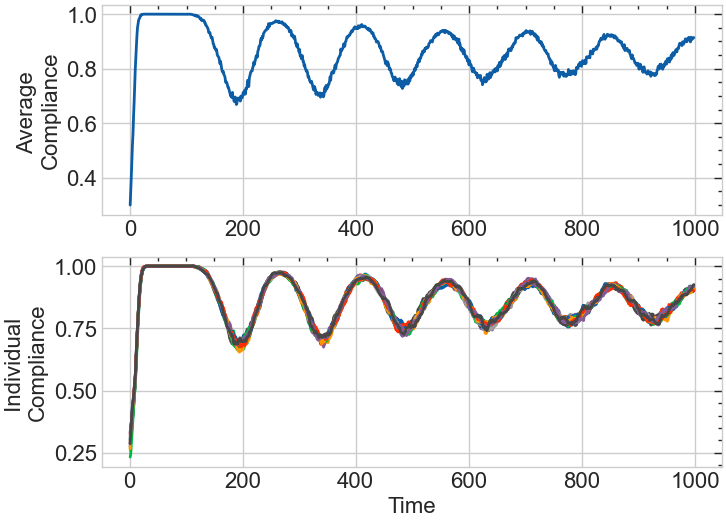}
\caption{Compliance control when the delay $m$ is a random variable sampled from a Gaussian distribution with parameters $\mu = 35$, $\sigma = 17.5$ and each agent has a chance of disconnecting from the network $\eta = 0.025$. Compliance is achieved and the use of personalised signals make it so that every agent contributes in a fair way.  }
\label{fig: wifi35}
\end{figure}

\begin{figure}
\includegraphics[width=1\columnwidth]{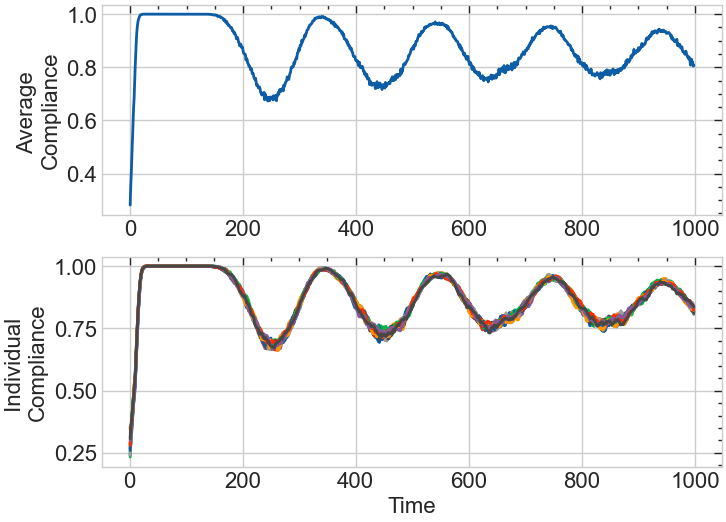}
\caption{Compliance control when the delay $m$ is a random variable sampled from a Gaussian distribution with parameters $\mu = 45$, $\sigma = 22.5$ and each agent has a chance of disconnecting from the network $\eta = 0.025$. Compliance is achieved and the use of personalised signals make it so that every agent contributes in a fair way.  }
\label{fig: wifi45}
\end{figure}

\section{Conclusions}
\label{sec: Conclusions and Future Research}

In this paper we explored the use of a feedback control system to regulate the behaviour of stochastic agents and to enforce the desired level of compliance, both globally and individually. The use of personalised feedback signals takes into account the behaviour of each agent and leads to fair regulation, with respect to each individual's base compliance $q_i$, whereas the global signal increases the robustness of the control system to malfunctions and malicious behaviour. We proved a theorem that establishes that the averaged compliance of each agent, under the proposed regulation scheme, will accumulate around the target compliance $Q^*$ and finally we validated our results through extensive Monte Carlo simulations. As per future lines of research we intend to provide theoretical results for the robustness of the proposed compliance control against malicious actors, explore different formulations of fairness to include, as an example, the economic status of each agent. Furthermore, we intend to extend our framework by using elements of game theory to take into account more complex scenarios.
\textcolor{black}{As a further strand of future work, we also wish to integrate our work into network simulators such as ns3. perhaps in combination with mobility simulators, to provide more detailed experimental validation of the proposed techniques. Finally, we recognise that the compliance work presented here involves the co-design of technology and behaviours and that enforcing compliance involves exploring and managing the appetite for risk in agents (including human decision makers). This suggests a strong connection to control strategies that involve simultaneous exploration and policy enforcement (such as reinforcement learning). We have already commenced work in this direction and future publications will report on this work.}

\section{Acknowledgements}

Ferraro, Zhao and Shorten are funded in part by the IOTA Foundation, and by Science Foundation Ireland grant 16/IA/4610 respectively.


\newpage
\appendix
\subsubsection{Proof of Lemma \ref{lem2}}\label{lem2-pf}
Both parts (a) and (b) of Lemma \ref{lem2} will follow by showing that $y(k)$ is well approximated by the solution of the 
following system of differential equations:
\be
\frac{d z_1}{d t} &=& h_1(z(t)) = w \, \left(p(z_{2}(t)) - z_1(t)\right) \label{ODE1} \\
\frac{d z_2}{d t} &=& h_2(z(t)) = \beta_0 \, w \, \left(Q^* - z_1(t) \right) \label{ODE2}
\ee
Defining $z(t) = (z_1(t), z_2(t))$ we write this system as
\be\label{ODE}
\frac{d z}{d t} = h(z(t)) 
\ee
It is easy to check that if $z_1(0) \in [0,1]$ then $z_1(t) \in [0,1]$ for all $t \ge 0$, and hence (since $\beta_0 \le 1$)
\be\label{bd-z}
\|\frac{d z}{d t}\| = \|h(z(t)) \| \le2 \, w
\ee
In order to make the connection between $y(k)$ and $z(t)$, we introduce
\be\label{def:why}
\wh{y}(t_k) = y(k) \quad \text{where $t_k = \epsilon \, k$, $k \ge -m$}
\ee
and then define $\wh{y}(t)$ for all $t \ge - m \epsilon$ using linear interpolation of the values at $\{t_k\}$.
We consider the solution of (\ref{recury2}) with some initial condition satisfying (\ref{norm-def2}):
\be\label{why1-1}
\wh{y}_1(t_k) &=& y_1(0) + \epsilon \, \sum_{j=0}^{k-1} h_1(y(j)) \nonumber \\
&=& y_1(0) +  \int_{0}^{t_k} h_1(y(\lfloor \epsilon^{-1} s\rfloor)) \, d s\nonumber \\
&=& y_1(0) + \int_{0}^{t_k} h_1(\wh{y}(s)) \, d s + \nonumber \\ 
 & & +\int_{0}^{t_k} \left[h_1(y(\lfloor \epsilon^{-1} s\rfloor)) - h_1(\wh{y}(s))\right] \, d s
\ee
For all $s \in [t_j, t_{j+1}]$, there is some $a \in [0,1]$ such that
\bee
\wh{y}(s) &=& a y(j) + (1-a) y(j+1) 
\eee
and therefore, using (\ref{bd-ydiff}),
\bee
&& \hskip-0.3in |h_1(y(\lfloor \epsilon^{-1} s\rfloor)) - h_1(\wh{y}(s))| \\
&& = |h_1(y(j)) - h_1(a y(j) +  (1-a) y(j+1))| \\
&& \le  2 \, w \, (1- a ) \, \| y(j+1) - y(j) \| \\
&& \le 4 \, w^2 \, \epsilon
\eee
Therefore we have the bound
\be\label{int-bd1}
| \int_{0}^{t_k} \left[h_1(y(\lfloor \epsilon^{-1} s\rfloor)) - h_1(\wh{y}(s))\right] \, d s |
\le 4 \, w^2 \, \epsilon \, t_k
\ee
By similar reasoning we find
\be\label{y2-bd8}
&& \hskip-0.5in \wh{y}_2(t_k) = y_2(0) + \int_{0}^{t_k} h_2(\wh{y}(s)) \, d s  \nonumber \\ 
 &&\hskip-0.2in  +\int_{m \epsilon}^{t_k} \left[h_2(y(\lfloor \epsilon^{-1} s\rfloor - m)) - h_2(\wh{y}(s - m \epsilon))\right] \, d s \nonumber \\ 
&&\hskip-0.2in - \int_{t_k-m \epsilon}^{t_k} h_2(\wh{y}(s)) \, d s + \epsilon \sum_{j-0}^{m-1} h_2(y(j-m))
\ee
and using (\ref{recury2}), (\ref{bd-ydiff}) and (\ref{norm-def2}) we get
\be\label{int-bd2}
&& \hskip-0.5in  | \int_{m \epsilon}^{t_k} \left[h_2(y(\lfloor \epsilon^{-1} s\rfloor - m)) - h_2(\wh{y}(s - m \epsilon))\right] \, d s | \nonumber \\ 
&& \hskip0.2in \le 2 \, w^2 \, \epsilon  \, t_k 
\ee
Now let $z(t)$ be the solution of (\ref{ODE}) with $z(0) = y(0)$.
For any $t \in [t_k,t_{k+1})$, using  (\ref{why1-1}), (\ref{y2-bd8}), (\ref{bd-ydiff}), (\ref{bd-z}), (\ref{int-bd1}) and (\ref{int-bd2}),
\be
&& \hskip-0.4in \| z(t) - \wh{y}(t)\| \nonumber \\
&&\hskip-0.4in \le  \| z(t_k) - \wh{y}(t_k)\| + \| z(t) -z(t_k) \| + \|\wh{y}(t) - \wh{y}(t_k)\| \nonumber \\
&& \hskip-0.4in \le  \int_{0}^{t_k} \, \|h(z(s))  - h(\wh{y}(s))\| \, d s \nonumber \\
&& + 6 \, w^2 \, \epsilon \, t_k + 2 \, w \, m \, \epsilon + 4 \, w \, \epsilon \nonumber \\
&& \hskip-0.4in \le  3 \,w \, \int_{0}^{t_k} \| z(s) - \wh{y}(s) \| d s + \epsilon \, \left(6 \, w^2 \, t_k + B_4\right) \nonumber \\
&& \hskip-0.4in \le 3 \,w \, \int_{0}^{t} \| z(s) - \wh{y}(s) \| d s + \epsilon \, \left(6 \, w^2 \, t + B_4\right)
\ee
where $B_4 = 2 \, w \, m + 4 \, w$.
The Gr\"onwall inequality now yields the bound
\be\label{diff1-bound}
\sup_{0 \le s \le t} \| z(s) - \wh{y}(s) \|  \le \epsilon \, \left(6 \, w^2 \, t + B_4\right) \, e^{3 \,w \, t}
\ee

\begin{lemma}\label{lem1}
\noindent For any $M \ge 0$, let $z(t)$ be the solution defined by
(\ref{ODE}) with initial conditions $0 \le z_1(0) \le 1$ and $|z_2(0)| \le M$.
\noindent{\bf a)} There is $\tau_1(M) < \infty$ such that 
\be\label{lem1:eq1}
0 \le z_2(t) \le 1 \,\, \text{some $t \le \tau_1(M)$}.
\ee
\noindent{\bf b)} 
There is $\tau_2(M) < \infty$ such that 
\be\label{lem1:eq2}
\| z(t) - (Q^*,Q^*) \| \le \frac{1}{2} \, \| z(0) - (Q^*,Q^*) \| \,\, \text{for all $t \ge \tau_2(M)$}.
\ee
\end{lemma}

Before proving Lemma \ref{lem1}, we use it to prove Lemma \ref{lem2}.
Given $M \ge 0$, define
\be
B_5 = \left(6 \, w^2 \, \tau+ B_4\right) \, e^{3 \,w \, \tau}
\ee
where $\tau = \tau_2(M)$ is the number defined in Lemma \ref{lem1} (b). For any
$\mathcal{Y}$ which is bounded at level $M$, let $y(k)$ be the solution of (\ref{recury2})
with initial condition $\mathcal{Y}$, and let $z(t)$ be the solution of (\ref{ODE}) with initial value $z(0) = y(0)$.
Then (\ref{diff1-bound}) implies that
\be
\|z(\tau) - \wh{y}(\tau) \|  \le B_5 \, \epsilon
\ee
Therefore
\be
&& \hskip-0.2in \|y(\lfloor\tau \epsilon^{-1}\rfloor) - (Q^*,Q^*) \| \nonumber \\
&& \hskip-0.3in \le  \| z(\tau) - (Q^*,Q^*) \| + \|z(\tau) - \wh{y}(\tau) \| + \|y(\lfloor\tau \epsilon^{-1}\rfloor) - \wh{y}(\tau)\| \nonumber \\
&& \hskip-0.3in \le  \frac{1}{2} \, \| z(0) - (Q^*,Q^*) \| + B_5 \, \epsilon + 2 \, w \, \epsilon \nonumber \\
&& \hskip-0.3in =  \frac{1}{2} \, \| y(0) - (Q^*,Q^*) \| + B_1 \, \epsilon
\ee
where $B_1 = B_5 + 2 \, w$, and this establishes Lemma \ref{lem2} (a).
For Lemma \ref{lem2} (b), note first that the bound $Y_1(k) \in [0,1]$ follows directly from (\ref{def-Y}).
Also the bound $|Y_2(k)| \le M_2$ will be implied by a bound on the duration of any excursion outside
the region $|Y_2| \le 2$. Indeed consider integers $u, v$ such that $|Y_2(u-1)| \le 2$ and
\be\label{def:v}
|Y_2(u+k)| > 2 \,\, \text{for $k=0,\dots,v$}
\ee
Let $y(k)$ be the solution of (\ref{recury2}) with initial conditions $Y(u-m:u)$,
and let $z(t)$ be the solution of (\ref{ODE}) with $z(0)=y(0)$.
It follows that $\xi_i(u+1)=0$, and so
$y_1(u+1)=Y_1(u+1)$ and $|Y_2(u+1)-y_2(u+1)| \le \epsilon^{3/2} \alpha_0$.
Also $\xi_i(u+k)=0$ for $k=1,\dots,v$, and so by similar reasoning it follows that
\be
y_1(u+k) &=& Y_1(u+k) \nonumber \\
|y_2(u+k) - Y_2(u+k)| & \le & \epsilon^{3/2} \alpha_0 \, k \label{y2-bd8}
\ee
for all $k \le v$ such that $\epsilon^{3/2} \alpha_0 \, k \le 1$.
Furthermore (\ref{lem1:eq1}) and (\ref{diff1-bound}) imply that
\be\label{y2-bd9}
|y_2(\lfloor \tau_1 \epsilon^{-1}\rfloor)| \le  1 + \epsilon \, B_6
\ee
where $B_6 = \left(6 \, w^2 \, \tau_1 + B_4\right) \, e^{3 \,w \, \tau_1}$ and
$\tau_1$ is defined in Lemma \ref{lem1} (a). For $\epsilon$ sufficiently small we
have $\epsilon^{3/2} \alpha_0 \,  \tau_1 \, \epsilon^{-1} \le 1$.
If we suppose that $v > \tau_1 \epsilon^{-1}$, then (\ref{y2-bd8}) and (\ref{y2-bd9}) would imply
that for $\epsilon$ sufficiently small
\be\label{Y2-bd7}
|Y_2(u + \lfloor \tau_1 \epsilon^{-1}\rfloor)| \le 1 + \epsilon^{1/2} \alpha_0 \, \tau_1 + \epsilon \, B_6 < 2
\ee
This contradicts (\ref{def:v}), so we conclude that $v \le \tau_1 \epsilon^{-1}$.
Therefore
\be
\sup_{0 \le k \le v} |Y_2(u+k)| \le 2 + \tau_1 \epsilon^{-1} \, \left(\epsilon \, w + \epsilon^{3/2} \, \alpha_0\right)
\ee
This holds for any excursion outside the region $|Y_2| \le 2$, and so we conclude that
\be
\sup_{k \ge 0} |Y_2(k)| \le 2 + \tau_1  \left(w + \epsilon^{1/2} \, \alpha_0\right)
\ee
which proves Lemma \ref{lem2} (b) with $M_2 = 2 + \tau_1  \left(w + \alpha_0\right)$.

\vskip0.2in
\subsubsection{Proof of Lemma \ref{lem1}}\label{lem1-pf}
Lemma \ref{lem1} (a) is trivial for $M \le 1$, so suppose $M > 1$.
In the region $\{z_2 \ge 1\}$ the solution of (\ref{ODE}) has the form
\be\label{ODE-sol2}
z_2(t) = z_2(0) - \beta_0 w t (1 - Q^*) + \beta_0 (1 - z_1(0))  (1 - e^{-w t})
\ee
Taking $z_2(0)=M$,  (\ref{ODE-sol2}) shows that $z_2(t_1) = 1$ for some $t_1 \le M \, \left(\beta_0 \,w \, (1 - Q^*)\right)^{-1}$,
and $z_2(t) \le M + \beta_0$ for all $0 \le t \le t_1$.
Similar reasoning for the region $\{z_2 \le -M\}$ with $z_2(0)=-M$ shows that 
$z_2(t_2) = 0$ for some $t_2 \le (M+1) \, \left(\beta_0 \,w \, Q^*\right)^{-1}$,
and $z_2(t) \ge - M - \beta_0$ for all $0 \le t \le t_2$.
This establishes (\ref{lem1:eq1}) with $\tau_1(M) = \max\{t_1,t_2\}$, and also shows that
\be\label{z-sup}
\sup_{t \ge 0} |z_2(t)| \le M + \beta_0
\ee

\medskip
To prove Lemma \ref{lem1} (b) we consider separately the cases $Q^* \ge 1/2$ and $Q^* < 1/2$.
Suppose first that $Q^* \ge 1/2$, and consider the function
\be
V(z) = \beta_0 \, \left(z_1 - Q^*\right)^2 + \left(z_2 - Q^*\right)^2
\ee
Note that $\dot{V} =  -2 \beta_0 w \left(z_1 - Q^*\right)^2$ in the square $[0,1]^2$, so
$V$ is a Lyapunov function in the square $[0,1]^2$. We define
\be
{\cal E} = \{z : V(z) \le \beta_0 \, \left(1 - Q^*\right)^2\}
\ee
The condition $\beta_0 \le 1$ implies that ${\cal E} \subset [0,1]^2$.
Therefore if the solution $z(t)$ enters ${\cal E}$ then it will remain thereafter inside the square $[0,1]^2$,
and thus its future evolution is determined by the linear system
\be\label{ODE1}
\dot{z_1} = w \, \left(z_{2} - z_1\right), \quad
\dot{z_2} = \beta_0 \, w \, \left(Q^* - z_1 \right)
\ee
This linear system can be written in matrix form as follows:
\be\label{ODE2}
\frac{d z}{ d t} = w\, A \, z + b, \, A =\bmx -1 & 1 \cr - \beta_0 & 0 \emx, \,
b = w \, \beta_0 \, Q^* \, \bmx 0 \cr 1 \emx.
\ee
It is easy to see that the matrix $A$ is stable with eigenvalues
\be
\lambda_{\pm} =  - \frac{1}{2} \pm \frac{1}{2} \sqrt{1 - 4 \beta_0^2}
\ee
Therefore the condition $\beta_0 > 0$ implies that $\lambda_{\pm}$ have negative real parts, and therefore
the solution of (\ref{ODE1})
converges exponentially to the fixed point $(Q^*,Q^*)$.
It remains to show that the solution $z(t)$ enters the set ${\cal E}$ within a bounded time.
Accordingly we define three closed subsets of the square $[0,1]^2$ as follows:
\bee
S_1 &=& \{z : 0 \le z_1 \le Q^*, \, 0 \le z_2 \le 1\} \\
S_2 &=& \{z : Q^* \le z_1 \le 1, \, Q^* \le z_2 \le 1\} \\
S_3 &=& \{z : Q^* \le z_1 \le 1, \, 0 \le z_2 \le Q^*\} 
\eee
It follows from the definition of ${\cal E}$ that
\bee
S_2 \cap S_3 =  \{z : Q^* \le z_1 \le 1, \, z_2 = Q^*\} \subset {\cal E}
\eee
Inspection of the system (\ref{ODE}) shows that the solution follows a trajectory that spirals
clockwise around the fixed point $(Q^*,Q^*)$. 
If $z(0) \in S_1$ then $z(t)$ must eventually reach $S_2$, either directly from $S_1$ or after an excursion into the region
$z_2 > 1$. We define $\tau_{12}$ to be the supremum over all starting points $z(0) \in S_1$ of the time until first
entering the set $S_2$. These times depend continuously on $z(0)$ and $S_1$ is compact, therefore $\tau_{12} < \infty$.
Similarly $\tau_{23} < \infty$ is the maximum time to reach $S_3$ starting from $S_2$, and $\tau_{31} < \infty$
is the maximum time to reach $S_1$ starting from $S_3$. 
Therefore there is $\tau_2' \le \tau_{12} + \tau_{23} + \tau_{31}$ such that starting from any point $z(0)$ in $[0,1]^2$,
the solution $z(t)$ will reach the interval $S_2 \cap S_3$ at some time $t \le \tau_2'$. 
If the trajectory starts at a point $z(0)$ satisfying $1 \le |z_2(0)| \le M$ then (\ref{lem1:eq1}) implies it must enter the square
$[0,1]^2$ 
before time $\tau_1$, and so must reach the interval $S_2 \cap S_3$ before time $\tau_2' + \tau_1$.

\medskip
Since $S_2 \cap S_3 \subset {\cal E}$
this shows that the solution $z(t)$ enters the set ${\cal E}$ within time $\tau_2'+ \tau_1$,
starting from any point in the region $|z_2| \le M$.
As noted before the system is a contraction in ${\cal E}$, so for any $r > 0$
there is some $\tau(r) < \infty$ such that for all $z(t) \in {\cal E}$,
\bee
\| z(t+s) - (Q^*,Q^*) \| \le r \, \| z(t) - (Q^*,Q^*) \| \, \text{for all $s \ge \tau(r)$}.
\eee
If $z(0) \in {\cal E}$ this establishes (\ref{lem1:eq2}) with $\tau_2(M)=\tau(1/2)$.
So assume that $z(0) \notin {\cal E}$. Let $R_1 = \{z : 0 \le z_1 \le 1, \, |z_2| \le M\}$
and define
\be
\rho = \sup_{z \in {\cal E}, w \in R_1 \setminus {\cal E}} \frac{\|z - y^*\|}{\|w - y^*\|}
\ee
Then for all $z(0) \in R_1\setminus {\cal E}$, we know that $z(t)$ reaches ${\cal E}$ by latest time $\tau_2'+ \tau_1$, and then
is contracted after time $\tau(r)$. So {for all $t \ge \tau(r) + \tau_2'+ \tau_1$ we have
\be
\| z(t) - (Q^*,Q^*) \| \le r \, \rho \, \| z(0) - (Q^*,Q^*) \|
\ee
Finally we choose $r \le (2 \rho)^{-1}$ and $\tau_2(M) = \tau(r) + \tau_2'+ \tau_1$ and conclude that
for all $z(0) \in R_1$,
\be\label{conv1-z}
\| z(t) - (Q^*,Q^*) \| \le \frac{1}{2} \, \| z(0) - (Q^*,Q^*) \| \, \text{for all $t \ge \tau_2(M)$}.\nonumber
\ee
A similar argument applies when $Q^* \le 1/2$, and this establishes Lemma \ref{lem1} (b).


\begin{thebibliography}{88}


\bibitem{accesspaper} Ferraro, P., King, C. and Shorten, B., "Distributed Ledger Technology for Smart Cities, the Sharing Economy and Social Compliance", {\em IEEE Access}, Vol. 6, pp. 62728 - 62746, 2018.

{\color{black}\bibitem{ML1} M. Loey, G. Manogaran, M. H. N. Taha, and N. E. M. Khalifa, A hybrid
deep transfer learning model with machine learning methods for face mask
detection in the era of the covid-19 pandemic {\em Measurement}, vol. 167, p.
108288, 2021.
\bibitem{ML2} E. Mbunge, S. Simelane, S. G. Fashoto, B. Akinnuwesi, and A. S. Metfula,
Application of deep learning and machine learning models to detect covid-19 face masks-a review {\em Sustainable Operations and Computers},
vol. 2, pp. 235245, 2021.
\bibitem{ML3} X. Kong, K. Wang, S. Wang, X. Wang, X. Jiang, Y. Guo, G. Shen, X. Chen and Q. Ni, "Real-time mask identification for covid-19: An edge-computing-based deep learning network", {\em IEEE Internet of Things Journal}, vol. 8, no. 21, pp. 15 92915 938, 2021.


\bibitem{ML4} A.Sanal and G. Udupa, "Machine learning based human body temperature measurement and mask detection by thermal imaging" {\em Third IEEE International Conference on Intelligent Computing Instrumentation and Control Technologies(ICICICT)},   vol 10, pp. 13391343, 2022.

\bibitem{Phan} Phan, T., Annaswamy, A. M., Yanakiev, D. and Tseng, E., ``A model-based dynamic toll pricing strategy for controlling highway traffic" {\em IEEE American Control Conference (ACC),} pp. 6245-6252, 2016.
{\color{black}\bibitem{Park} B. Park, K. Amasyali, Y. Chen and M. Olama, "Hierarchical Transactive Control of Flexible Building Loads Under Distribution LMP,"  {\em 2022 IEEE Power and Energy Society Innovative Smart Grid Technologies Conference (ISGT)}, pp. 1-5, 2022.
\bibitem{Nudell} Nudell, T. R., Brignone, M., Robba, M., Bonfiglio, A., Ferro, G., Delfino, F. and Annaswamy, A. M., "Distributed control for polygeneration microgrids: A Dynamic Market Mechanism approach", {\em Control Engineering Practice,} Volume 121, 2022.}

\bibitem{Soylemezgiller} Soylemezgiller, F., Kuscu, M. and Kilinc, D., ``A traffic congestion avoidance algorithm with dynamic road pricing for smart cities. In Personal Indoor and Mobile Radio Communications ", {\em IEEE 24th Annual International Symposium on Personal, Indoor, and Mobile Radio Communications (PIMRC) }, pp. 2571-2575, 2013.

\bibitem{Bui} Bui, K. T., Huynh, V. A., and Frazzoli, E., `` Dynamic traffic congestion pricing mechanism with User-Centric considerations", {\em 15th International IEEE Conference on Intelligent Transportation Systems}, pp. 147-154, 2012.
\bibitem{Annaswamy}Annaswamy, A. M., Guan, Y., Tseng, H. E., Zhou, H., Phan, T. and Yanakiev, D, ``Transactive Control in Smart Cities",  {\em Proceedings of the IEEE}, Vol. 106, No. 4, pp. 518-537, 2018.
\bibitem{Widergren} Widergren, S., Fuller, J., Marinovici, C. and Somani, A., ``Residential transactive control demonstration", {\em IEEE PES Innovative Smart Grid Technologies Conference (ISGT)}, pp. 1-5, 2014.
\bibitem{Huang} Huang, P., Kalagnanam, J., Natarajan, R., Hammerstrom, D., Melton, R., Sharma, M. and Ambrosio, R., ``Analytics and transactive control design for the pacific northwest smart grid demonstration project"  {\em First IEEE International Conference on Smart Grid Communications (SmartGridComm),} pp. 449-454, 2010.
\bibitem{Kotb} Kotb, A. O., Shen, Y. C., Zhu, X. and Huang, Y., ``iParker-A New Smart Car-Parking System Based on Dynamic Resource Allocation and Pricing",  {\em IEEE Transactions on Intelligent Transportation Systems,} Vol. 17, No. 9, pp. 2637-2647, 2016.
\bibitem{Yao}Yao, Y. and Zhang, P., Transactive control of air conditioning loads for mitigating microgrid tie-line power fluctuations", {\em IEEE PES General Meeting}, pp. 1-1,2016
\bibitem{Katipamula} Katipamula, S., ``Smart buildings can help smart grid: Transactive controls"{\em IEEE PES Innovative Smart Grid Technologies (ISGT)}, pp. 1-1, 2012.
\bibitem{Junjie} Junjie, H. U., Guangya, Y. A. N. G., Koen, K. O. K., Yusheng, X. U. E. and Bindner, H. W., ``Transactive control: a framework for operating power systems characterized by high penetration of distributed energy resources"  {\em Journal of Modern Power Systems and Clean Energy}, Vol. 5, No. 3, pp. 451-464, 2017.
\bibitem{Li} Li, J., Lin, X., Nazarian, S. and Pedram, M., ``CTS2M: concurrent task scheduling and storage management for residential energy consumers under dynamic energy pricing"  {\em IET Cyber-Physical Systems: Theory and Applications}, Vol. 2, No. 3, pp. 111-117, 2017.

\bibitem{Chekired} Chekired, D. A., Khoukhi, L. and Mouftah, H. T., ``Decentralized cloud-SDN architecture in smart grid: A dynamic pricing model"  {\em IEEE Transactions on Industrial Informatics}, Vol. 14, No. 3, pp. 1220-1231, 2018.
\bibitem{Bejestani} Bejestani, A. K., Annaswamy, A. and Samad, T., ``A hierarchical transactive control architecture for renewables integration in smart grids: Analytical modeling and stability", {\em IEEE Transactions on Smart Grid}, Vol. 5, No. 4, pp. 2054-2065, 2014.
\bibitem{Hao} Hao, H., Corbin, C. D., Kalsi, K. and Pratt, R. G., ``Transactive control of commercial buildings for demand response"  {\em IEEE Transactions on Power Systems}, Vol. 32, No. 1, pp. 774-783, 2017.

\bibitem{Popov} Popov, S., ``The Tangle-Version 1.4.3", available at \url{https://iota.org/IOTA_Whitepaper.pdf}, April 2018.
{\color{black}\bibitem{TangleNew} S. Muller, A. Penzkofer, N. Polyanskii, J. Theis, W. Sanders and H. Moog, "Tangle 2.0 Leaderless Nakamoto Consensus on the Heaviest DAG," in {\em IEEE Access}, Vol. 10, pp. 105807-105842, 2022,}
\bibitem{Tangle} Popov, S., Saa, O. and Finardi, P., ``Equilibria in the Tangle" {\em arXiv preprint} arXiv:1712.05385, 2017.

\bibitem{IoTPaper1} Cullen, A., Ferraro, P., King, C. and Shorten, R., ''On the resilience of dag-based distributed ledgers in iot applications'', {\em IEEE Internet of Things Journal}, Vol. 7, No. 8, pp. 7112-7122, 2020.
\bibitem{IoTPaper2} Moschella, M., Ferraro, P., Crisostomi, E. and Shorten, R., ''Decentralized Assignment of Electric Vehicles at Charging Stations Based on Personalized Cost Functions and Distributed Ledger Technologies'', {\em IEEE Internet of Things Journal}, 2021.
\bibitem{TaCPaper} Ferraro, P., King, C. and Shorten, R., ''On the stability of unverified transactions in a DAG-based Distributed Ledger'', {\em IEEE Transactions on Automatic Control}, Vol. 65, No.9, pp. 3772-3783, 2019.
\bibitem{ethereum1} Buterin, V., ''A next-generation smart contract and decentralized application platform'', {\em white paper}, Vol. 3 No.37, 2014.
{\color{black}\bibitem{ethereum2} Wang, Z., Jin, H., Dai, W., Choo, K. K. R., and Zou, D., "Ethereum smart contract security research: survey and future research opportunities", {\em Frontiers of Computer Science}, Vol. 15, No. 2, pp. 1-18 2021.}
\bibitem{Rob-Mon} Robbins, H., and Monro, S., ``A stochastic approximation method'', {\em Ann. Math. Statistics}, 22:400D407, 1951.
\bibitem{Borkar} Borkar, V.~S., ``Stochastic approximation a dynamical systems viewpoint'', Hindustan Book Agency, New Delhi, 2008.
\bibitem{ns3} Riley, G. F. and Henderson, T. R.,``The ns-3 Network Simulator" ,  in {\em Modeling and Tools for Network Simulation}, editors - Wehrle, K. and Günes, M. and Gross, J., Spinger, pp. 15-34, 2010.

}
\end{thebibliography}
\end{document}